\documentclass[preprint,11pt,authoryear,nonatbib]{article}

\usepackage[utf8]{inputenc}
\usepackage[backend=biber,style=apa]{biblatex}
\usepackage{setspace}
\usepackage{graphicx} 
\usepackage{amsmath, amssymb}
\usepackage{mathtools}
\usepackage[margin=1in]{geometry}
\usepackage[toc,page]{appendix}
\usepackage{pgfplotstable}
\usepackage{booktabs}
\usepackage{siunitx}
\usepackage{multirow}
\usepackage{lscape}
\usepackage{pdflscape}
\usepackage{calc} 
\usepackage{enumitem} 
\usepackage[hidelinks]{hyperref} 
\usepackage{float}
\usepackage{cleveref}
\usepackage{lscape}
\usepackage{tabularx}
\usepackage{ragged2e}
\usepackage{authblk}
\usepackage{makecell}
\usepackage{standalone}

\usepackage{amsthm}
\usepackage{xspace} 
\usepackage{caption}
\usepackage{subcaption}
\usepackage[para]{threeparttable}
\usepackage{changes}

\usepackage{accents}

\usepackage{lscape}
\usepackage{xcolor}
\usepackage{color} 
\usepackage{colortbl}

\usepackage[linesnumbered,ruled,vlined,noresetcount]{algorithm2e}
\RestyleAlgo{ruled}
\SetKwComment{Comment}{/* }{ */}
\usepackage{booktabs}
\usepackage[normalem]{ulem}
\usepackage{ulem}
\usepackage{cancel} 
\allowdisplaybreaks

\usepackage{subfiles} 
\usepackage{tikz}
\usetikzlibrary{shapes, arrows, positioning, arrows.meta, shapes.geometric, fit, backgrounds}

\tikzstyle{process} = [rectangle, draw, text centered, minimum height=2em]
\tikzstyle{connector} = [draw, -latex]
\tikzstyle{terminator} = [rectangle, draw, text centered, rounded corners, minimum height=2em]

\tikzset{
  block/.style args={#1}{
    rectangle, rounded corners,
    draw=#1!70!black,
    fill=#1!12,
    thick,
    align=center,
    minimum width=3.7cm,
    minimum height=0.85cm,
    font=\small
  },
  smallblock/.style args={#1}{
    rectangle, rounded corners,
    draw=#1!70!black,
    fill=#1!10,
    thick,
    align=center,
    minimum width=3.2cm,
    minimum height=0.75cm,
    font=\scriptsize
  },
  decision/.style={
    diamond,
    aspect=2,
    draw=green!50!black,
    fill=green!10,
    thick,
    align=center,
    inner sep=1pt,
    font=\scriptsize
  },
  arrow/.style={-{Latex[length=2mm]}, thick},
  dashedarrow/.style={-{Latex[length=2mm]}, thick, dashed},
  groupbox/.style args={#1}{
    draw=#1!70!gray,
    fill=#1!5,
    rounded corners,
    dashed,
    thick,
    inner xsep=8pt,
    inner ysep=8pt
  }
}
\usepackage{pgfplots} 
\pgfplotsset{compat=1.18}
\usepgfplotslibrary{groupplots}

\newcommand{\linePool}{\ensuremath{\mathcal{L}}}
\newcommand{\frequencySet}{\ensuremath{\mathcal{F}}}

\newcommand{\qosheuristic}{SAMCF-heuristic\xspace}

\newcommand{\stopSet}{\ensuremath{N}\xspace}
\newcommand{\stopSetCG}{\ensuremath{\mathcal{N}}\xspace}

\newcommand{\stopSetCGA}{\ensuremath{\mathcal{N}^{A}}\xspace}
\newcommand{\stopSetCGE}{\ensuremath{\mathcal{N}^{E}}\xspace}
\newcommand{\stopSetCGT}{\ensuremath{\mathcal{N}^{T}}\xspace}
\newcommand{\stopSetCGL}{\ensuremath{\mathcal{N}^{L}}\xspace}
\newcommand{\stopn}{\ensuremath{n}\xspace}
\newcommand{\linePT}{\ensuremath{l}\xspace}

\newcommand{\ODset}{\ensuremath{\mathcal{D}}\xspace}
\newcommand{\ODor}{\ensuremath{s_d}\xspace}
\newcommand{\ODde}{\ensuremath{t_d}\xspace}
\newcommand{\ODw}{\ensuremath{w_d}\xspace}
\newcommand{\od}{\ensuremath{d}\xspace}

\newcommand{\ModeSet}{\ensuremath{M}\xspace}
\newcommand{\mode}{\ensuremath{m}\xspace}

\newcommand{\Pset}{\ensuremath{\mathcal{P}}\xspace}
\newcommand{\PsetOD}{\ensuremath{\mathcal{P}_d}\xspace}
\newcommand{\PsetODPT}{\ensuremath{\PsetOD^{PT}}\xspace}
\newcommand{\PsetODedge}{\ensuremath{\mathcal{P}_d(a)}\xspace}
\newcommand{\PSetLineFreq}{\ensuremath{\mathcal{P}_d(l,f)}\xspace}

\newcommand{\arc}{\ensuremath{a}\xspace}

\newcommand{\arcSet}{\ensuremath{\mathcal{A}}\xspace}
\newcommand{\arcSetDrive}{\ensuremath{\mathcal{A}^{IVT}}\xspace}
\newcommand{\arcSetTransfer}{\ensuremath{\mathcal{A}^{T}}\xspace}
\newcommand{\arcSetAcces}{\ensuremath{\mathcal{A}^{A}}\xspace}
\newcommand{\arcSetEgress}{\ensuremath{\mathcal{A}^{E}}\xspace}
\newcommand{\arcSetOptOut}{\ensuremath{\mathcal{A}^{O}}\xspace}
\newcommand{\linePoolSol}{\ensuremath{\mathcal{L}^{*}}\xspace}

\newcommand{\ticketFare}{\ensuremath{R^{\text{fare}}}\xspace}
\newcommand{\subsidy}{\ensuremath{R^{\text{sub}}}\xspace}
\newcommand{\incentive}{\ensuremath{R}\xspace}

\newcommand{\probODPT}{\ensuremath{p_\od^{PT}}\xspace}
\newcommand{\utilODALT}{\ensuremath{u^{ALT}_\od}\xspace}
\newcommand{\utilODPT}{\ensuremath{u^{PT}_\od}\xspace}

\newcommand{\curSol}{\ensuremath{s}\xspace}
\newcommand{\newSol}{\ensuremath{s^{\prime}}\xspace}
\newcommand{\bestSol}{\ensuremath{s^{\star}}\xspace}
\newcommand{\initSol}{\ensuremath{s^{0}}\xspace}

\newcommand{\fcurSol}{\ensuremath{J}\xspace}
\newcommand{\fnewSol}{\ensuremath{J^{\prime}}\xspace}
\newcommand{\fbestSol}{\ensuremath{J^{\star}}\xspace}

\newcommand{\destroySet}{\ensuremath{\Theta^{-}}\xspace}
\newcommand{\repairSet}{\ensuremath{\Theta^{+}}\xspace}
\newcommand{\repairOperator}{\ensuremath{n^{+}}\xspace}
\newcommand{\destroyOperator}{\ensuremath{n^{-}}\xspace}

\newcommand{\probRepair}[1]{\ensuremath{p^{+}_{#1}}\xspace}
\newcommand{\probDestroy}[1]{\ensuremath{p^{-}_{#1}}\xspace}
\newcommand{\reward}[1]{\ensuremath{\sigma^{#1}}\xspace}
\newcommand{\score}[1]{\ensuremath{\psi^{#1}}\xspace}
\newcommand{\attempts}[1]{\ensuremath{\nu^{#1}}\xspace}
\newcommand{\weightDestroy}[1]{\ensuremath{\rho^{-}_{#1}}\xspace}
\newcommand{\weightRepair}[1]{\ensuremath{\rho^{+}_{#1}}\xspace}

\newcommand{\reactionFactor}{\ensuremath{\gamma}\xspace}
\newcommand{\coolingFactor}{\ensuremath{\delta}\xspace}

\newcommand{\baselineShare}{\ensuremath{p_d^{\mathrm{obs}}}\xspace}
\newcommand{\cng}{CNG\xspace}

\usepackage[all]{nowidow} 
\usepackage[protrusion=true,expansion=true]{microtype} 

\title{An ALNS Heuristic for Large-Scale Line Planning with Mode Choice and Line Generation}

\author[1]{Siv Marie Cartland Hansen}
\author[1]{Richard Martin Lusby}
\author[2]{Rowan Hoogervorst}
\author[1]{Otto Anker Nielsen}

\affil[1]{Department of Technology, Management and Economics, Technical University of Denmark, Akademivej 358, Kongens Lyngby, 2800, Denmark}
\affil[2]{Department of Strategy, Entrepreneurship and Operations, EDHEC Business School, 24 avenue Gustave Delory, 59057 Roubaix, France}

\date{}

\begin{document}

\maketitle

\begin{abstract}
Demand responsiveness is an important consideration in public transport line planning, as network design and service quality influence passenger demand. However, accounting for this interaction further complicates an already challenging combinatorial optimization problem. To address this challenge, we propose a scalable Adaptive Large Neighborhood Search (ALNS) algorithm for large-scale line planning with endogenous demand. The algorithm jointly optimizes lines and frequencies while accounting for passenger mode choice, passenger assignment, and vehicle capacities. Candidate lines are generated dynamically throughout the search, and solutions are evaluated using an embedded evaluation procedure for passenger assignment and demand estimation, together with a dedicated local search procedure for frequency optimization.
 The proposed methodology is evaluated on the public transport network of Odense, Denmark, comprising approximately 1,800 origin-destination pairs. Computational results demonstrate the applicability of the approach to realistic, large-scale instances. The optimized networks concentrate resources on fewer, higher-frequency services, reducing average headways from approximately 41 minutes to 6.8–-13 minutes while substantially increasing public transport ridership. Furthermore, the results show that network design is highly sensitive to assumptions regarding passenger behavior, highlighting the importance of carefully calibrated demand models when incorporating demand responsiveness into line planning.
\end{abstract}

\section{Introduction}

Public transport (PT) planners routinely face the challenge of designing networks that provide attractive and efficient passenger services. Network design decisions directly affect passenger accessibility, travel times, and the overall attractiveness of public transport. As transport authorities seek to improve mobility and support a transition towards more sustainable travel behavior, understanding how network design influences travel demand has become increasingly important. Consequently, transit network design methods should account for the interaction between service provision and traveler behavior, recognizing that changes in service quality influence passengers' mode choices and, ultimately, public transport ridership.

Line planning is a central component of transit network design, involving the selection of which lines to operate and at what frequencies. The resulting set of operated lines and their associated frequencies is commonly referred to as the \textit{line concept}. Most line planning models assume exogenous demand and select lines from a predefined set of candidate lines. Yet, even under these simplifying assumptions, exact solution methods struggle to scale to instances of realistic size within practical computation times. In reality, public transport demand depends on the quality of the service offered, as passengers may change their mode choice in response to variations in travel times, waiting times, and network coverage. Accounting for such behavioral responses substantially increases computational complexity by coupling line planning decisions with passenger assignment and demand estimation through a mode choice model. Moreover, restricting the optimization to a predefined line pool may limit the quality of the resulting network design. Allowing lines to be generated dynamically eliminates the need to enumerate candidate lines in advance but considerably increases the solution space. Consequently, there remains a need for efficient solution approaches that account for demand responses to service quality while remaining applicable to realistic networks. 

In this paper, we present a scalable matheuristic for line planning with demand determined through a logit-based mode choice model, in which passengers choose between public transport and a single competing mode. The proposed solution method is based on an Adaptive Large Neighborhood Search (ALNS) algorithm that integrates mixed-integer programming and problem-specific heuristic operators to address the computational challenges arising from endogenous demand, passenger assignment, and dynamic line generation. Adopting a system-optimal perspective, the algorithm seeks network designs that improve overall system performance while accounting for passengers' responses to service changes. The methodology is evaluated on the regional public transport network surrounding the Danish city of Odense, comprising approximately 1,800 origin-destination (OD) pairs and representing a problem size that is beyond the capabilities of exact optimization approaches. 

The paper makes three main contributions. First, we develop a competitive and scalable algorithm for line planning that accounts for endogenous demand and dynamic line generation. Second, we provide insights from a large-scale real-world case study, illustrating how accounting for mode choice influences network design and passenger outcomes. Finally, we evaluate the passenger path assignment procedure embedded in the framework and assess its performance relative to more detailed route choice models. Our numerical results indicate that the optimized networks concentrate resources on fewer, high-frequency lines, reducing average headways from around 41 minutes under current operations to 6.8--13 minutes. When induced demand is taken into account, the model selects to invest heavily in network intensification, leading to higher operating costs but also capturing a much larger market share, from approximately 8\% of total demand to up to 32\% as the public transport system becomes more competitive with outside options. However, the resulting designs are sensitive to the behavioral assumptions, highlighting the importance of careful calibration of the demand response. 

The paper begins by positioning the work within the existing literature in Section~\ref{sec:lit}. Section~\ref{sec:problem} then introduces the line planning problem and its mathematical formulation. Building on this, Section~\ref{sec:ALNS} describes the proposed ALNS algorithm, detailing the destroy and repair operators together with the local search procedure used to adjust service frequencies. The framework is subsequently applied to the Danish case study in Section~\ref{sec:results}, where its computational performance and the resulting planning insights are examined on realistic-sized instances. Finally, Section~\ref{sec:conclusion} summarizes the main findings and discusses directions for future research.

\section{Literature review}
\label{sec:lit}

	As one of the core steps within the typical public transport planning cycle, line planning has been extensively studied in the literature.
	We refer to the literature review by \textcite{schobel_line_2012} and the more recent preprint by \textcites{schmidt_planning_2024} for a general introduction.
	As discussed in these literature reviews, two distinguishing factors between line planning approaches are the way in which passenger behavior is modeled and the freedom available in choosing lines and frequencies.
	We discuss these factors further in the remainder of this section.
	We note that the line planning problem is also closely related to the network design problem, which in its most common form comprises the generation and selection of lines, while line planning traditionally encompasses line selection and frequency setting for a predefined line pool.
	For an overview on the network design literature, we refer to \textcite{ceder_bus_1986}, \textcite{guihaire_transit_2008}, and \textcite{duran-micco_survey_2022}. 
	
	\subsection{The modeling of passenger behavior}
	A first distinction between line planning approaches lies in how passenger behavior is modeled.
	Considering the difficulty of line planning problems, initial approaches focused on finding lines with minimum operating costs considering a fixed assignment of passengers to edges or paths.
	For example, by finding lines such that a given minimum (and possibly maximum) frequency is achieved on each edge in the network, or by first assigning passengers onto a shortest path through the public transport network (before line choice) and minimizing the resulting travel times.
	Alternatively, others have looked at choosing lines such to give the largest amount of passengers a direct link without transfers from their origin to destination \autocite{bussieck_optimal_1997,suman2019improvement}.
	
	To better capture the effect of the chosen lines on passenger travel times, one needs to take into account the effect that the selected lines have on the paths that passengers can take within the network.
	It is common to use the Change \& Go network as introduced by \textcite{schobel_line_2006} in this context, which encodes the travel options of passengers in terms of traveling on lines and transferring between them.
	Many approaches then assume that the passengers' paths are chosen within the model according to a system optimum, i.e., passengers follow a route that leads to the lowest overall system cost even if this might mean a (small) increase in travel time for some of the passengers.
	Examples include the papers of \textcite{schobel_line_2006}, \textcite{borndorfer_column-generation_2007}, and \textcite{hurk_shuttle_2016}.
	Note that a system optimum assignment fits well into a traditional mixed-integer programming approach and that the system optimum assignment can, in some cases, be close to a free route choice by passengers \parencite[see, e.g., the results in][]{hansenExactAlgorithmPublic2026}.
	
	Some recent approaches have further integrated route choice for passengers, i.e., passenger paths can differ from the shortest or the system optimal path.
	These papers often use a bi-level optimization framework, where the lower level ensures passengers choose a route corresponding to their shortest path or based on a logit-based distribution \parencite[see, e.g.,][]{cancela2015mathematical, goerigk_line_2017, cipriani_bus_2012}.
	Here,  \textcite{goerigk_line_2017} linearize the lower-level formulation through either dualization or adding implicit shortest path constraints, while the linearization by \textcite{hartleb_modeling_2023} is based on the assumption that passengers distribute evenly over routes with similar journey time.
	The latter two approaches are, however, applied to relatively smaller railway instances.
	
	A second modeling choice related to passenger behavior lies in whether public transport demand is considered fixed or dependent on the offered level of service.
	The majority of papers assume that a public transport modal split is determined beforehand, i.e., that the total public transport demand is not (significantly) affected by the offered line concept.
	See, e.g., \textcite{schobel_line_2006,borndorfer_column-generation_2007,hurk_shuttle_2016}.
	More recently, papers have started to treat passenger mode choice as endogenous in the model, capturing that high service quality can induce demand and that low service quality can instead mean passengers switch to other transport modes.
	Here, \textcite{hansenExactAlgorithmPublic2026} consider a maximum travel time above which passengers are assumed to use a transport mode other than public transport. 
	\textcite{bertsimas_data-driven_2021} assume that the number of passengers increases linearly with the best available travel time for that passenger group.
	Moreover, \textcite{hartleb_modeling_2023} model a logit-based mode choice between rail and an alternative mode and show that considering mode choice is important for both operator profit and attracting more passengers in their numerical experiments.
	
	\subsection{Selection of lines and frequencies}
	Line planning papers also differ in terms of the level of choice and detail considered in determining the public transport network.
	First, in terms of the possible lines that can be chosen.
	Most commonly, line planning papers assume a line pool to be given, implying a sequential process in which a set of possible lines is first generated in a line generation step and this line pool is then input to line planning.
	While working with a line pool gives fine-grained control of the possible lines, generating a good set of lines is a challenging problem by itself, and the choice of line pool can have a significant impact on the quality of the downstream line concept found in line planning \autocite{gattermann_line_2017}.
	In contrast to using a predefined line pool, some papers have looked at generating promising lines within the line planning stage.
	This is sometimes referred to as line planning with full pool.
	For example, \textcites{borndorfer_column-generation_2007} and \textcite{bertsimas_data-driven_2021} use a column generation approach where lines are generated dynamically, while \textcites{canca_integrated_2019} adapt lines as part of a neighborhood-based heuristic.
	
	Related to the possibilities in choosing lines are the possible frequencies that can be chosen.
	Some approaches assume that the frequency can be chosen as a continuous number, allowing any frequency between a given lower and upper bound.
	See, e.g., the papers by \textcites{borndorfer_column-generation_2007} and \textcites{lu_line_2025}.
	While a continuous frequency is easier to integrate into exact approaches, frequencies are often picked in practice from a set of easy-to-remember options.
	For example, many bus lines operate every 10, 15, 20, or 30 minutes.
	Such a discrete frequency has been taken into account by, e.g., \textcites{bertsimas_data-driven_2021} and \textcites{hurk_shuttle_2016}.
	To better handle the model complexity that comes with handling a larger number of discrete frequencies,  \textcite{hansenExactAlgorithmPublic2026} have proposed an exact iterative solution approach that dynamically expands the number of considered frequencies.
	
	\subsection{Heuristics for line planning}
	The line planning literature is split into both exact and heuristic approaches.
	The computational boundary of problems that can still be solved exactly has been shifting over time and depends strongly on the studied line planning problem, e.g., in terms of the demand modeling considerations discussed above.
	Recent exact approaches are, though, able to solve (reasonably) large real-life instances.
	For example, in a recent preprint, \textcites{lu_line_2025} solve line planning instances with consideration of crowding effects, using a diving heuristic, that contain up to 277 stations.
	\textcites{bertsimas_data-driven_2021} solve bus network instances with up to 233 stops and 2042 edges with some consideration of mode choice.
	To solve such large instances, exact approaches often make some simplifying assumptions.
	For example, \textcites{bertsimas_data-driven_2021} assume public transport demand scales linearly with travel time instead of assuming a full logit-based mode choice.
	Altogether, this suggests there is still a trade-off between realism in modeling and the magnitude of instances that can be solved by exact approaches.
	
	To solve large-scale real-life instances, and especially ones with more detailed passenger consideration, heuristics have thus also been extensively used for line planning.
	Different metaheuristic frameworks have been explored within this context, with population-based metaheuristics being popular.
	For example, \textcite{bielli2002genetic} study fitness functions for the use within a genetic algorithm for bus network design, and \textcites{cipriani_bus_2012} use a genetic algorithm to solve a bus network design problem with elastic demand.
	\textcites{szeto_hybrid_2012} propose a hybrid enhanced artificial bee colony algorithm for bus network design.
	\textcite{fan2006simulated} propose a simulated annealing heuristic for transit network design.
	A memetic algorithm, i.e., a genetic algorithm that integrates with a solver, is developed by \textcites{duran-micco_designing_2022} to solve a real-life case study for the city of Utrecht in the Netherlands.
	
	Another stream of heuristic approaches focuses on neighborhood-based metaheuristic frameworks.
	Of particular interest, \textcites{canca_adaptive_2017} propose an ALNS heuristic for line planning with full pool and mode choice in a railway setting.
	Some of their neighborhoods, such as extending existing lines in the solution, align with the neighborhoods we use in our ALNS heuristic.
	They extend their work in \textcite{canca_integrated_2019}, where a distribution of passengers over paths is considered instead of an all-or-nothing assignment to a single path.
	A variable neighborhood search heuristic is proposed by \textcites{iliopoulou2022variable} for transit network design (without frequency setting).
	Moreover, a tabu search algorithm for network design is proposed by \textcites{fan2008tabu}, which iterates between the generation of new transit lines and the assignment of frequencies and passengers to routes.

    Numerous heuristics for line planning also fall in the category of matheuristics, i.e., heuristics that include a mathematical programming component.
	These include the earlier mentioned memetic algorithm of \textcite{duran-micco_designing_2022} and ALNS algorithm of \textcite{canca_integrated_2019}.
	Moreover, \textcite{gatt_solving_2025} propose a matheuristic based on column generation and enumeration and test it for a case study for the bus network of a mid-size French city.
	Also our algorithm falls into this category of matheuristics, as we use column generation for passenger assignment within the evaluation of a line concept.
	
	\subsection{Summary}
    This paper contributes to the line planning literature by proposing a highly scalable ALNS algorithm for a realistic version of the line planning problem with both mode choice and dynamic line generation.
	As shown by the schematic overview of the line planning literature given in \autoref{tab:literaturereview}, only  \textcites{canca_integrated_2019} and \textcites{bertsimas_data-driven_2021} have so far looked at combining mode choice and dynamic line generation as we do in this paper.
	Compared to \textcites{bertsimas_data-driven_2021}, we model mode choice in a more realistic way by assuming that mode choice is made according to a logit distribution instead of assuming that demand scales linearly with weighted travel time.
	Moreover, while \textcites{canca_integrated_2019} consider very similar line planning properties and also use an ALNS algorithm, their passenger assignment step requires considerable computation time, limiting the size of instances that can be solved.
	Instead, we show that our paper can solve large-scale, realistic bus networks, like our case study for the bus network within the city of Odense.
		
	\begin{table}[htbp]
				
		\caption{An overview of the line planning literature and the positioning of our ALNS-based approach. Abbreviations: SO: System Optimum, RC: Route Choice; MC: Mode Choice; LP: Line Pool; LG: Line Generation, Co.: Continuous, Di.: Discrete, Ex.: Exact Approach; Heur.: Heuristic.}
        \label{tab:literaturereview}
		\centering
		
		\definecolor{Gray}{gray}{0.94}
		\newcolumntype{a}{>{\columncolor{Gray}}c}
		
		\begin{threeparttable}
			\begin{tabular}{lcacacacacac}
				\toprule
				& \multicolumn{3}{c}{Demand} & \multicolumn{2}{c}{Lines} & \multicolumn{2}{c}{Freq.\ } & \multicolumn{2}{c}{Algorithm} & \multicolumn{2}{c}{Domain} \\
				\cmidrule(rl){2-4} \cmidrule(rl){5-6} \cmidrule(rl){7-8} \cmidrule(rl){9-10} \cmidrule(rl){11-12}
				Paper & SO & RC & MC & LP & LG & Co.\ & Di.\ & Ex.\ & Heur.\ & Rail & Bus \\
				\midrule
				\textcites{schobel_line_2006} & \checkmark & & & \checkmark & & & \checkmark &  \checkmark  & & \checkmark & \\
				\textcite{borndorfer_column-generation_2007} & \checkmark & & & & \checkmark & \checkmark & & \checkmark & & & \\
				\textcite{cipriani_bus_2012} &  & \checkmark & \checkmark & \checkmark  &  & \checkmark &   &  & \checkmark &   & \checkmark \\
				\textcite{hurk_shuttle_2016} & \checkmark & & & \checkmark & & &  \checkmark &  \checkmark & & &  (\checkmark)\tnote{a} \\
				\textcite{goerigk_line_2017} & & \checkmark & & \checkmark & & & \checkmark \tnote{b} & \checkmark & (\checkmark) & \checkmark & \\
				\textcite{canca_adaptive_2017} & \checkmark & & \checkmark & & \checkmark & & \checkmark & & \checkmark \tnote{e} & \checkmark & \\
				\textcite{canca_integrated_2019} & \checkmark & & \checkmark & & \checkmark & & \checkmark & & \checkmark \tnote{e} & \checkmark & \\
				\textcite{bertsimas_data-driven_2021} & \checkmark & & \checkmark \tnote{c} & & \checkmark & & \checkmark & \checkmark & & & \checkmark \\
				\textcite{hartleb_modeling_2023} & & \checkmark & \checkmark & \checkmark & & \checkmark & & \checkmark & & \checkmark & \\
				\textcite{lu_line_2025} & \checkmark & & & \checkmark & & \checkmark & &  \checkmark & & \checkmark & \\
				\textcites{gatt_solving_2025} & \checkmark\tnote{d} & & & & \checkmark & & \checkmark \tnote{b}  & & \checkmark \tnote{e} & & \checkmark \\
				\textcite{hansenExactAlgorithmPublic2026} & \checkmark & & (\checkmark)\tnote{f}  & \checkmark & & & \checkmark & \checkmark & & & \checkmark \\
				This paper & \checkmark & & \checkmark & & \checkmark & & \checkmark & & \checkmark \tnote{e} & & \checkmark \\
				\bottomrule
			\end{tabular}
			\begin{tablenotes}
				\item[a] Rail-replacement minibuses
				\item[b] Any integer frequency (between bounds) is allowed
				\item[c] Linear scaling of demand to travel time
				\item[d] Upper bound on path deviations to shortest path
				\item[e] Matheuristic
				\item[f] OD pairs will choose other transport mode if travel time exceeds threshold
			\end{tablenotes}
	\end{threeparttable}
	\end{table}

\section{Problem description and model}
\label{sec:problem}

We consider the problem of designing a bus network within a multimodal PT system, where PT competes with an alternative mode such as a private car, cycling, walking, or another non-PT option. While passengers may use other public transport services, we focus on optimizing the bus network and treat other PT modes as fixed. For brevity, we refer to all boarding locations --- whether bus stops or rail stations --- uniformly as stops. We assume that the underlying network, including the set of candidate stops and connections, is fixed and given. In addition, operational costs, vehicle capacities, travel demand patterns, and travel times for the alternative mode are assumed to be known for the planning period considered (typically a peak hour). 

The problem is to select a set of bus lines and their operating frequencies that optimize a system-level objective that accounts for passenger generalized travel cost and operator costs, while capturing the endogenous effect of service quality on the amount of passenger demand that the PT network can attract and accommodate. The notation used in this section is collected in \autoref{tab:notation}. 

\begin{table}[htbp!]
    \centering
    \caption{Notation and terminology}
        \label{tab:notation}
    \renewcommand{\arraystretch}{0.97}
    \begin{tabular}{ll}
        \toprule
        \textbf{Sets} & Description \\
        \midrule
        \ModeSet & set of PT modes \\
        \stopSet & set of candidate stops \\
        $E$ & set of bidirectional edges \\
        $E_\mode \subseteq E$ & set of edges for mode $\mode \in \ModeSet$\\
        $\linePool$ & set of candidate lines \\
        $E_\linePT \subseteq E$ & set of edges traversed by line $\linePT \in \linePool$ \\
        $\frequencySet$ & set of all candidate frequencies \\
        $\frequencySet_\linePT \subseteq \frequencySet$ & set of candidate frequencies for line $\linePT \in \linePool$ \\
        \ODset & set of OD pairs \\
        $\arcSet$ & set of directed arcs in the \cng \\
        $\arcSetDrive \subseteq \arcSet$ & arcs representing in-vehicle travel \\
        $\arcSetDrive_\linePT \subseteq \arcSetDrive$ & in-vehicle arcs associated with line $\linePT$ \\
        $\arcSetAcces \subseteq \arcSet$ & boarding arcs \\
        $\arcSetEgress \subseteq \arcSet$ & egress arcs \\
        $\arcSetTransfer \subseteq \arcSet$ & transfer arcs \\
        $\arcSetOptOut \subseteq \arcSet$ & alternative mode arcs \\
        \PsetOD & set of paths for OD pair $\od \in \ODset$ \\
        $\PsetOD(a) \subseteq \PsetOD$ & set of paths of OD pair \od traversing arc $a \in \arcSet$ \\
        \PsetODPT  $\subseteq \PsetOD$ & set of paths of OD pair \od using only PT arcs \\
        $ \PSetLineFreq \subseteq  \PsetOD$ & Set of paths  of OD pair \od using line $l \in \linePool$ at frequency $f \in  \frequencySet_l$ \\
        \midrule
        \textbf{Parameters} & Description \\
        \midrule
        \ODw & demand of OD pair $\od \in \ODset$ (passengers/period) \\
        $\ODor \in \stopSet$ & origin stop of OD pair $\od$ \\
        $\ODde \in \stopSet$ & destination stop of OD pair $\od$ \\
        $c_\arc$ & generalized cost of arc $\arc \in \arcSet$ \\
        \utilODALT & generalized cost of the alternative mode for OD pair $\od$ \\
        $\alpha_{\od},\, \beta_{\od}$ & logit model parameters for OD pair $\od \in \ODset$\\
        $\phi_{\linePT f}$ & minimum vehicles to operate line $\linePT$ at frequency $f$ \\
        $\delta_{\mode}$ & passenger capacity per vehicle of mode $\mode \in \ModeSet$ \\
        $h_{\linePT}$ & fixed cost of operating line $\linePT  \in \linePool$ \\
        $c_m$ & cost per vehicle of mode $\mode \in \ModeSet$ \\
        \incentive & revenue collected per passenger, including fare \ticketFare and subsidy \subsidy \\ 
        $B$ & operational budget \\
        \midrule
        \textbf{Variables} & Description \\
        \midrule
        $y_{\linePT f}$ & whether to operate line $\linePT \in \linePool$ at frequency  $f \in \frequencySet_\linePT$ \\
        $z_{\linePT}$ & number of vehicles assigned to line $\linePT$ \\
        $x_{pd}$ & the fraction of passengers of OD pair $d$ on path $p\in \PsetOD$ \\
        \utilODPT & expected generalized cost of PT for OD pair \od \\
        \probODPT & fraction of demand of OD pair $\od$ willing to use PT \\
        \bottomrule
    \end{tabular}
\end{table}

\subsection{Public transport network and line representation} 
The Public Transport Network (PTN), i.e., the underlying public transport infrastructure network, is represented by a multigraph $G = (\stopSet, E) = (\stopSet, \bigcup_{\mode \in \ModeSet} E_\mode)$.
Here, \stopSet is a set of candidate stops, and each $E_\mode$ is a set of bidirectional edges that can be served by mode $\mode \in \ModeSet$. 
Each edge $e \in E$ connects a pair of stops and is associated with a travel time $t_e$. We refer to the subgraph of $G$ which can be served by mode $\mode$ as $G_\mode = (\stopSet, E_\mode)$. 

A line $\linePT$ of mode $m$ is defined as a simple path through $G_m$, representing a bidirectional service between two terminal stops without repeated stops. 
Let $\mode(\linePT)$ return the mode of line \linePT. 
All feasible lines are assumed to be given in a set $\linePool = \bigcup_{\mode \in \ModeSet} \linePool_\mode$ (a line pool). We denote by $E_\linePT \subseteq E$ the edges traversed by line \linePT. Each line $\linePT \in \linePool$ can be operated at one of a discrete set of candidate frequencies $\frequencySet_\linePT$.
The frequency specifies the number of vehicle departures per unit time and can equivalently be expressed as the headway of a line which gives the number of minutes between consecutive vehicles. 
A solution specifies an operated line set $\linePoolSol \subseteq \linePool$, where each line $\linePT \in \linePoolSol$ is assigned exactly one frequency $f^* \in \frequencySet_\linePT$. 

\subsection{Operating costs and capacities}
Operating costs are modeled primarily through fleet size, which is taken as the main driver of variable operating expenses. Each vehicle of mode $\mode \in \ModeSet$ incurs a unit cost $c_\mode$, representing the flexible cost of operating a vehicle over the planning period, including costs such as driver wages, energy consumption, and routine maintenance. In addition, a fixed cost $h_{l}$ is incurred when operating a line $l\in\linePool$, capturing line-specific overhead. 

The cost and capacity of a line are jointly determined by its length and the selected operating frequency. 
For each line, we define a round-trip travel time that captures the time needed to serve the line in both directions, including terminal turnaround time. The required fleet size, $\phi_{lf}$, denotes the number of vehicles needed to operate line \linePT at frequency $f\in \frequencySet_l$, and is increasing in the service frequency.
The capacity of line $l$ is then determined by the vehicle capacity $\delta_{\mode(l)}$ and the number of vehicles $\phi_{lf}$ assigned to the line during the planning period. 
In the model, capacity constraints are enforced at the arc and line level, ensuring that passenger flows on each arc associated with a line do not exceed the line capacity. 

Passenger transport also generates revenue for the operating agency. 
We assume a fixed per-passenger incentive, \incentive, comprising both a per-passenger fare revenue \ticketFare and subsidy \subsidy. 
This can be extended to OD-specific fares, distance-based pricing, or alternative subsidy schemes.

\subsection{Passenger routing}
To determine the level of service provided by a line concept, it is essential to determine the routing of passengers through the network.
Here, we assume the interest of passengers to travel is given by OD pairs $\ODset$.
Each $\od \in \ODset$ denotes an OD pair between an origin $\ODor\in N$ and destination stop $\ODde \in N$. 
For each pair \od, the number of passengers who wish to travel is given by $w_d$. 
This quantity is defined for the specific planning period, usually a peak hour. 
While total travel demand is given and fixed, the share of demand willing to use public transport depends on the service offered relative to an alternative mode, which requires determining both the service level in the public transport system and the mode choice through a demand function.

\subsubsection{Passenger flows and the Change \& Go-network}
Passenger flows are modeled using the Change \& Go network, which was first introduced by \textcite{schobel_line_2006}.
This network corresponds to a directed graph \cng = (\stopSetCG, \arcSet) that allows for representing boarding a vehicle, driving on vehicles, and transfers between PT services. 
Following \textcite{hansenExactAlgorithmPublic2026}, we consider a frequency-expanded version of the \cng network and further extend the graph to capture demand that cannot be routed within the PTN.

Nodes in the \cng are linked to stops and the node set is defined as $\stopSetCG:= \stopSetCGA  \cup \stopSetCGE  \cup  \stopSetCGT \cup \stopSetCGL $.
Here, $\stopSetCGA$ and $\stopSetCGE$ contain access and egress nodes for each $n\in \stopSet$ that is also an origin or destination stop for at least one OD pair in $\ODset$, \stopSetCGT contains transfer nodes at stops $n\in \stopSet$ served by at least two lines, and \stopSetCGL contains line-specific stop nodes, one for each stop along the route of a line.

Arcs then represent activities taken by passengers, where the arc set is defined as $\arcSet: = \arcSetDrive \cup \arcSetAcces \cup \arcSetEgress \cup \arcSetTransfer\cup \arcSetOptOut$.
In-vehicle arcs \arcSetDrive, that represent passengers using a line on a segment, are defined for each $e \in E_l$ of each line $\linePT \in \linePool$, one for each direction. 
Arcs specific to line $l \in \linePool$ are denoted by  $\arcSetDrive_l \subseteq \arcSetDrive$. 
Access arcs \arcSetAcces connect origin nodes in \stopSetCGA to relevant line-stop nodes in \stopSetCGL and represent passengers boarding a vehicle from their origin. 
Passengers access a specific line at a specific frequency, and an arc is therefore defined for each frequency of the line. 
Egress arcs \arcSetEgress, representing passengers arriving at their destination, connect a line-stop node in \stopSetCGL to an egress node in \stopSetCGE for every stop $\stopn\in \stopSet$, and every line $\linePT \in \linePool$ that stops at this stop, and every possible frequency $f \in \frequencySet_{l}$.
Transfer arcs \arcSetTransfer represent passengers transferring at a stop. Therefore, we define arcs from every node in \stopSetCGL to the transfer node of that stop in \stopSetCGT, a separate arc for each frequency. Opt-out arcs \arcSetOptOut represent travelers choosing an alternative mode of transport. These arcs connect nodes in \stopSetCGA to nodes in \stopSetCGE with one direct arc for each unique OD pair $(\ODor, \ODde)\in \ODset$. 

Passenger flows are defined as path flows in the \cng. Let \PsetOD denote the set of feasible paths for OD pair \od. For convenience, we define a PT-specific path set as $\PsetODPT = \{\, p \in \Pset_\od \mid p \cap \arcSetOptOut = \emptyset \,\}$,  the paths using arc $a \in \arcSet$ as $\Pset(a)$, and the paths using some line $l \in \linePool$ at frequency $f \in \frequencySet_l$ as $\PSetLineFreq$. Each path $p$ is associated with a generalized travel cost, dependent on the arcs used
\begin{equation}
    c_p = \sum_{\arc \in p} c_a   \qquad \forall p \in \Pset,
\end{equation}
where $c_a$ is the cost of an arc $a \in \arcSet$, which accounts for different components of journey time for the PT arcs or directly as $\utilODALT$ for $a \in \arcSetOptOut$. Furthermore, path costs for all paths $p \in \PsetODPT$ also include the ticket fare $\ticketFare$, which can be added as a constant or reflected in either the cost of access or egress arcs. 

Routing decisions are determined centrally in a system-optimal manner and therefore do not represent individual route choice behavior or user-equilibrium conditions. Although the objective seeks to minimize passenger inconvenience, demand associated with a given OD pair may still be distributed across multiple --- potentially non-shortest --- paths when capacity constraints are binding and the marginal reduction in passenger cost does not justify the additional capacity required to serve all flow.

\subsubsection{Passenger demand}
Public transport is assumed to compete with a single aggregated alternative mode, whose generalized cost is taken as fixed and independent of the public transport line concept. 
Although this abstraction ignores interactions between modes (e.g., shared roadway capacity), it is commonly adopted in strategic planning models \parencite[see, e.g.,][]{hartleb_modeling_2023}.
We model the proportion of travelers associated with OD pair $\od$ choosing public transport using the logit-function often used in transport applications \parencite[]{ben1985discrete}:
\begin{equation}
    \probODPT \leq \frac{1}{1 + e^{\left(\alpha_{\od} - \beta_{\od}\bigl(\utilODALT - \utilODPT \bigr)\right)}} \qquad \forall \od \in \ODset,
    \label{eq:logit}
\end{equation}
where \utilODALT is the fixed generalized cost of the alternative mode, $\alpha_{\od}\geq 0$ and $\beta_{\od}\geq 0$ are logit parameters, and \utilODPT is the experienced generalized cost of public transport for OD pair $\od \in \ODset$. 
We define the experienced PT service level for an OD as the weighted average generalized cost: 
\begin{equation}
      \utilODPT = \frac{\sum_{p \in \PsetODPT}c_p x_{pd}}{\sum_{p \in \PsetODPT}x_{pd}}  \qquad \forall \od \in \ODset,  
\end{equation}
where $x_{pd}$ is the flow on path $p$ for OD $\od \in \ODset$. 
The logit function links service quality, defined as the experienced generalized cost, to demand: the routing determines the level of service, which in turn influences the number of passengers willing to use public transport.

In particular, while the demand function determines the maximum share of travelers willing to use public transport, hard capacity constraints may prevent all willing demand from being accommodated. This modeling choice is motivated partly by computational tractability in the integrated routing and demand estimation problem and partly by the strategic planning perspective adopted in this study. We take the position that at this level of planning, operators should be able to design the most efficient network consistent with the objective function, without having to introduce additional capacity.

\subsection{Model formulation}\label{sec:model}
Using the above notation, we formulate a mixed‑integer nonlinear program (MINLP) extending \textcite{hansenExactAlgorithmPublic2026} and \textcite{hurk_shuttle_2016}. 
In particular, the interaction between flow assignment and capacity decisions follows \textcite{hansenExactAlgorithmPublic2026}, but we now model the maximum number of passengers willing to choose PT using \eqref{eq:logit} rather than a travel-time threshold. Furthermore, we replace the weighing parameter $\lambda$ with a PT-specific parameter \subsidy, which applies only to PT paths rather than to all alternatives.
The model jointly selects a line concept and a passenger routing, considering demand responsiveness and capacity constraints. 
To model the selection of opening lines, binary variables $y_{lf}$ determine whether a line $l \in \linePool$ is opened at frequency $f \in \frequencySet_l$.
The number of vehicles used on a line $l \in \linePool$ is then determined by the variables $z_l$.
Continuous variables are used to represent the earlier introduced flow of passengers $x_{pd}$, experienced PT service level $\utilODPT$, and proportion of travelers willing to use public transport $\probODPT$. 

  \begin{align}
    \min \;
    \left(
    \sum_{\od\in \ODset} 
    \sum_{p \in \PsetOD}
    c_p \ODw x_{p\od}\right)  +
    \left( 
    \sum_{l \in \linePool} 
    c_{\mode(l)} z_{l} + 
    \sum_{l \in \linePool} 
    \sum_{f \in \frequencySet_l} 
    h_{l} y_{lf} - 
     \sum_{\od\in \ODset} 
    \sum_{p \in \PsetODPT}
    \incentive \ODw x_{p\od} 
    \right) \label{obj} 
     \end{align}
     \vspace{-0.5cm}
subject to \vspace{-0.5cm}
     \begin{align}
     \sum_{p \in \PsetOD} x_{p\od} &= 1 && \forall \od \in \ODset \label{con1}  \\
        \sum_{\od\in \ODset} \sum_{p \in \PsetODedge} \ODw x_{p\od} &\leq \delta_{\mode(l)} z_{l}   && \forall l \in \linePool, \forall a \in \arcSetDrive_l  \label{con2}  \\
         \sum_{p \in \PSetLineFreq} x_{p\od} &\leq y_{lf}  && \forall \od \in \ODset, \forall l \in \linePool, \forall f \in \frequencySet_l \label{con4} \\
       \sum_{f \in \frequencySet_l} y_{lf} &\leq 1 && \forall l \in \linePool{} \label{con5} \\
       z_l &\geq y_{lf} \phi_{lf} && \forall l \in \linePool, \forall f \in \frequencySet_l \label{con6}\\
      \utilODPT &= \frac{\sum_{p \in \PsetODPT}c_p x_{pd}}{\sum_{p \in \PsetODPT}x_{pd}} && \forall \od \in \ODset  \label{conPTService}\\
   \probODPT &\leq \frac{1}{1 + e^{\!\left(\alpha_{\od} - \beta_{\od}\bigl(\utilODALT - \utilODPT \bigr)\right)}} && \forall \od \in \ODset  \label{conPTShare} \\
        \sum_{p\in \PsetODPT} x_{pd} &\leq \probODPT && \forall \od \in \ODset  \label{conPTFlow} \\
      \sum_{l\in \linePool}\left( c_{\mode(l)} z_{l} + \sum_{f\in\mathcal{F}_l} h_l y_{lf}\right) &\leq B  && \label{conBudget}    \\
       x_{p\od}&\geq 0 && \forall \od \in \ODset, \forall p \in \PsetOD \label{domain_x}\\   
       z_l &\geq 0 && \forall l \in \linePool \label{domain_z} \\  
       y_{lf} &\in \{0,1\} && \forall l \in \linePool, \forall f \in \frequencySet_l \label{domain_y}  \\
              \utilODPT &\geq 0 && \forall \od \in \ODset \label{domain_u}\\ 
                     \probODPT &\in [0,1] && \forall \od \in \ODset \label{domain_p}
\end{align}

The objective \eqref{obj} represents a trade-off between service and costs. 
The first term minimizes the total generalized travel cost experienced by travelers, while the second term accounts for operating costs and revenue. 
Constraints \eqref{con1} ensure that all demand associated with each OD pair is fully assigned to feasible paths in the \cng, possibly including paths representing the use of the alternative mode.  
Capacity constraints \eqref{con2} restrict passenger flows on each (directed) arc of a line to not exceed the capacity of the line. 
Constraints \eqref{con4} link passenger flows to line–frequency decisions by allowing passengers to use a line only if the corresponding frequency is activated. 
Constraints \eqref{con5} enforce that at most one frequency is selected for each line, with no frequency being selected meaning that the line is not opened, while constraints \eqref{con6} ensure that sufficient vehicles are assigned to operate a line at the chosen frequency.
Constraints \eqref{conPTService} define the experienced generalized PT cost for each OD pair as the average travel cost. 
This service level directly affects demand through the logit‑based constraints \eqref{conPTShare}, which determine an upper bound on the share of passengers willing to use public transport.
Constraints \eqref{conPTFlow} ensure that the demand routed on PT does not exceed this bound. Finally, constraint \eqref{conBudget} imposes an upper limit on total operating resources, and constraints \eqref{domain_x}–\eqref{domain_p} specify the domains of the decision variables.

The model contains two sources of nonlinearity, arising from the definitions of PT utility \eqref{conPTService} and PT demand \eqref{conPTShare}. 
Furthermore, although the formulation assumes that both the set of all feasible lines and the set of all candidate passenger paths are given, explicitly generating these sets in full is computationally intractable for large-scale instances. 
As a result, solving the problem to optimality (or even finding a feasible solution) is only possible for small-scale instances, which motivates the heuristic solution method presented in the following section.

\section{Solution method} \label{sec:ALNS}

We propose an ALNS algorithm to solve the model. 
ALNS is a metaheuristic belonging to the class of Very Large Neighborhood Search (VLNS) algorithms defined in \textcite{ahujaSurveyVeryLargescale2002}, which are based on the idea that searching large neighborhoods leads to high-quality local optima. 
ALNS particularly builds on the idea introduced in the Large Neighborhood Search (LNS) algorithm \autocite{shawUsingConstraintProgramming1998} to search large neighborhoods by destroying and repairing solutions using a set of operators, and complements it with a learning mechanism that guides the selection of operators based on their historical performance.
The method was first introduced in \textcite{ropkeAdaptiveLargeNeighborhood2006} for variants of the vehicle routing problem, and has since proven to be a competitive framework for a wide range of (transportation) problems; we refer to \textcites{pisinger_large_2019} for an overview. 
	
\subsection{Overview of the ALNS algorithm}
In our implementation, we employ four destroy operators in combination with three repair operators. 
Following the approach of \textcite{canca_integrated_2019}, we decouple the problem into two subproblems and solve them iteratively. The network design and operational decisions are handled directly by the ALNS operators, whereas demand estimation and passenger routing are managed through a dedicated evaluation subroutine. 
We note that while enumerating the complete set of feasible transit lines is computationally intractable for large-scale instances, we assume the availability of a predefined line pool that contains a set of reasonable lines. This line pool is used by one operator to add pre-screened lines, however, the method still allows to generate new unseen lines beyond the line pool.

The destroy operators modify the current solution by shortening or removing lines, while the repair operators construct new lines or extend existing ones. 
At each iteration, operators are selected probabilistically based on their recent performance, allowing the algorithm to adaptively favor more effective strategies over time.
The algorithm maintains three solutions: the current solution \curSol, the temporary candidate solution \newSol, and the incumbent solution \bestSol. Their objective values are denoted \fcurSol, \fnewSol, and \fbestSol, respectively. The evaluation of a candidate solution requires solving the passenger assignment problem and estimating mode choice simultaneously. This is accomplished using the proposed heuristic for the Service-Aware Multi-Commodity Flow (SAMCF) problem. As introduced in the preprint by \textcite{hansenColumnGenerationbasedFixedpoint2026}, the SAMCF problem jointly determines the demand to be served and its routing, where attracted demand is elastic and governed by a logit choice model, while routing is subject to hard capacity constraints.

\begin{figure}[t]
\centering

\begin{subfigure}[t]{0.48\textwidth}
\centering
  \resizebox{1.0\textwidth}{!}{%
\begin{tikzpicture}

\node[block=white] (start) {Start};

\node[smallblock=white, fill = none, below=of start] (init)
{\textbf{Initialize}\\
\scriptsize initial solution \initSol, best solution \bestSol, \\ operator scores};

\node[smallblock=white, fill = none, below=1.0cm of init] (destroy)
{Select and apply destroy operator\\\scriptsize  $\Rightarrow$ partial solution};

\node[smallblock=white,fill = none, below=of destroy] (repair)
{Select and apply repair operator\\
\scriptsize $\Rightarrow$ candidate solution \newSol};

\node[smallblock=white, below=of repair, xshift=-2cm] (localsearch)
{\textbf{Local search}\\
\scriptsize (improvement phase)};

\node[smallblock=gray, right=0.8cm of localsearch] (eval)
{\textbf{Evaluation step}\\
\scriptsize (\qosheuristic)};

\node[smallblock=white,fill = none, xshift=2cm, below=of localsearch] (accept)
{\textbf{Acceptance criterion}\\
\scriptsize accept or reject \newSol};

\node[smallblock=white,fill = none, below=of accept] (update)
{\textbf{Update}
\scriptsize \curSol if accepted, \\ \bestSol if improved, and operator scores};

\node[block=white, below=1.0cm of update] (stop)
{\textbf{Return}\\
\scriptsize best solution \bestSol};

\draw[arrow] (start) -- (init);
\draw[arrow] (init) -- (destroy);
\draw[arrow] (destroy) -- (repair);
\draw[arrow] (repair.south) -- ++(0,-0.35) -| (localsearch.north);
\draw[arrow] (localsearch.east) -- 
node[above,font=\scriptsize] {}
(eval.west);
\draw[arrow] (eval.west) -- 
node[above,font=\scriptsize] {}
(localsearch.east);
\draw[arrow] (localsearch.south) -- ++(0,-0.5) -| 
node[pos=0.25,left,font=\scriptsize] {}
(accept.north);
\draw[arrow] (accept) -- (update);
\draw[arrow] (update) -- (stop);

\draw[arrow]
  (update.west) -- ++(-1.5,0) |- node[pos=0.015,right,font=\scriptsize] {iterate}
  (destroy.west);

\begin{scope}[on background layer]
\node[groupbox=gray, fill = none, inner xsep=20pt, inner ysep=15pt, fit=(destroy)(repair)(localsearch)(eval)(accept)(update)] (alnsbox) {};
\end{scope}

\node[anchor=north west, font=\scriptsize\bfseries, text=gray!50!black]
at (alnsbox.north west) {ALNS main loop};
\end{tikzpicture}}
\caption{ALNS procedure}
\label{fig:alns-overview}
\end{subfigure}%
\hfill
\begin{subfigure}[t]{0.48\textwidth}
\centering
\resizebox{1.0\textwidth}{!}{
\begin{tikzpicture}

\node[block=gray, fill=none, minimum width=3.2cm] (input)
{Start};
 
\node[smallblock=white, fill=none, below=of input, minimum width=5.0cm] (init_qos)
{Initialize RMP \scriptsize (Set up MCF problem \\ with demand estimate $w_\od \ \forall \ \od\in\ODset$)};

\node[smallblock=white, fill=none, below=of init_qos, minimum width=5.0cm] (rmp)
{Solve RMP\\
\scriptsize obtain objective $z^{MCF}$ and duals $\pi$};

\node[smallblock=white, fill=none, below=of rmp, minimum width=5.0cm] (pricing)
{Solve pricing problem\\
\scriptsize (search for improving columns)};

\node[decision, below=of pricing, xshift=2cm, fill=none] (improve)
{Improving\\column?};

\node[smallblock=white, fill=none, left=1.0cm of improve, minimum width=3.3cm] (addcol)
{Add column\\to RMP};

\node[smallblock=white, fill=none, below=1.0cm of improve, xshift=-2cm, minimum width=5.0cm] (repair)
{\scriptsize Evaluate demand wrt. $\mathbf{x}$ and demand function\\
\scriptsize $\Rightarrow$ passenger costs $z^{REP}$};

\node[smallblock=white, fill=none, below=1.0cm of repair, xshift=-1cm, minimum width=3.0cm] (update_demand)
{Update demand \\ estimate $w_\od \ \forall \ \od\in\ODset$};

\node[decision, right=0.5cm of update_demand, fill=none] (stop)
{Converged?};

\node[block=white, fill=none, below=of update_demand, xshift=1cm,minimum width=3.2cm] (return)
{Return $z^{REP}$ and routing};

\draw[arrow] (input) -- (init_qos);
\draw[arrow] (init_qos) -- (rmp);
\draw[arrow] (rmp) -- (pricing);
\draw[arrow] (pricing.south) -- ++(0,-0.5) -| (improve.north);

\draw[arrow] (improve.west) -- node[above, font=\scriptsize] {yes} (addcol.east);
\draw[arrow] (improve.south) -- ++(0,-0.5) -| 
node[pos=0.5,left,font=\scriptsize] {no}
(repair.north);

\draw[arrow] ([xshift=-40pt]addcol.north) |- (rmp.west);

\draw[arrow] (repair.south) -- ++(0,-0.5) -|
(stop.north);
\draw[arrow] (stop.south) -- ++(0,-0.5) -| node[pos=0.5,left,font=\scriptsize] {yes}
(return.north);
\draw[arrow] (stop.west) -- node[above, font=\scriptsize] {no} (update_demand.east);

\draw[arrow]
  (update_demand.west) -- ++(-1.5,0) |- node[pos=0.015,right,font=\scriptsize] {iterate}
  (init_qos.west);
  
\begin{scope}[on background layer]
\node[groupbox=gray, solid, inner xsep=30pt,  fit=(input)(init_qos)(rmp)(pricing)(improve)(addcol)(return)] (cgbox) {};
\end{scope}

\node[anchor=north west, font=\scriptsize\bfseries, text=gray!40!black]
at (cgbox.north west) {Evaluation step};

\begin{scope}[on background layer]
\node[groupbox=gray, inner ysep = 20pt, fit=(rmp)(pricing)(improve)(addcol)] (cgbox) {};
\end{scope}

\node[anchor=north west,  font=\scriptsize\bfseries, text=gray!40!black]
at (cgbox.north west) {Column generation};
\end{tikzpicture}
}
\caption{Evaluation procedure}
\label{fig:evaluation-cg}
\end{subfigure}
\caption{Schematic overview of the ALNS algorithm and the embedded evaluation procedure.}
\label{fig:alns-and-evaluation}
\end{figure}
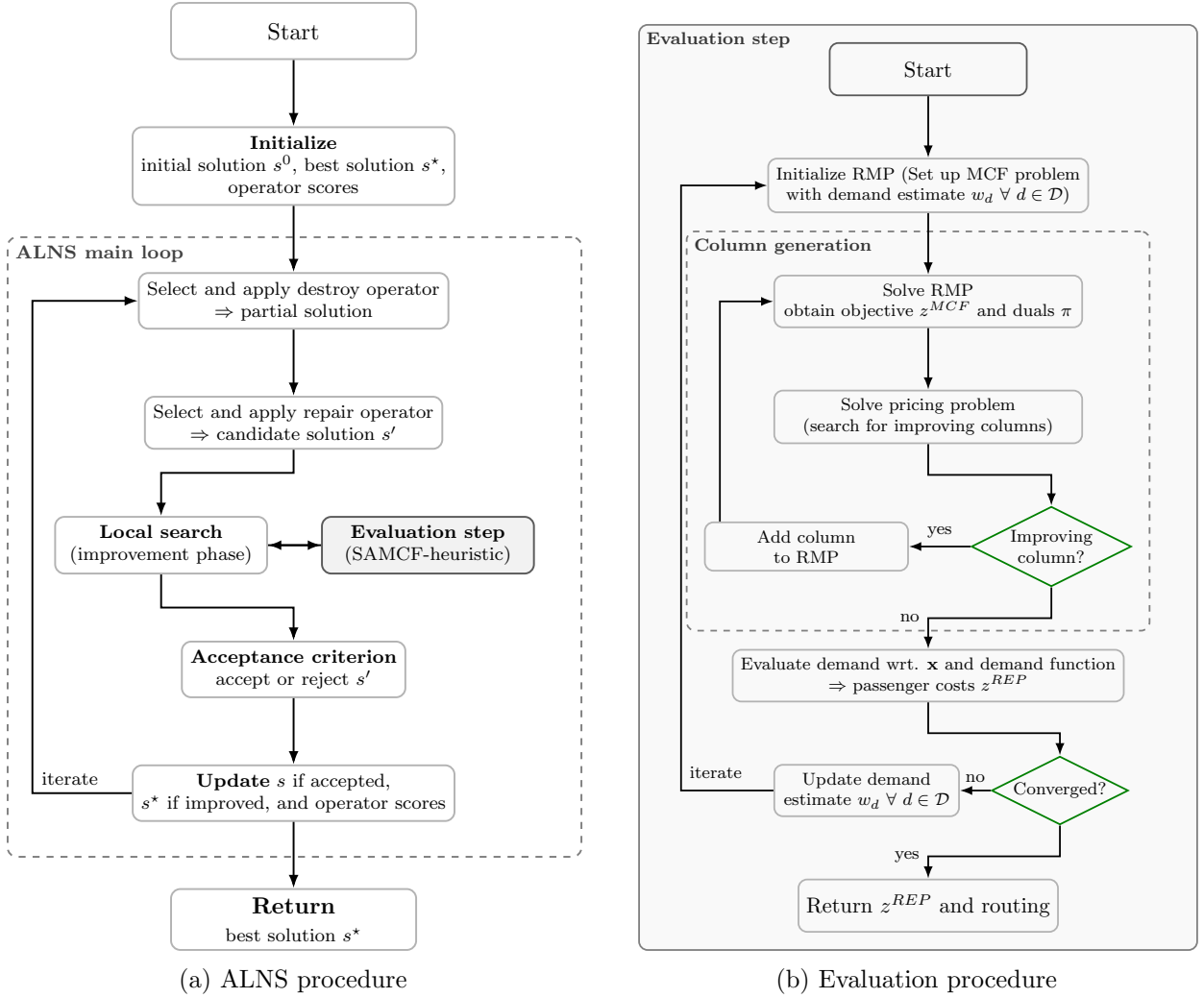

The evaluation step repeatedly solves a linear program, and one repair operator additionally solves a small Mixed-Integer Linear Program (MILP), meaning that the overall approach can also be classified as a matheuristic. 
The overall structure of our algorithm, together with an overview of the evaluation procedure, is presented in \autoref{fig:alns-and-evaluation}. 

In the following, we describe the components of the ALNS algorithm. \autoref{sec:neighborhoods} introduces the neighborhood operators, while \autoref{sec:localsearch} presents the local search procedure for adjusting line frequencies. The adaptive update of operator selection probabilities is described in \autoref{sec:adaptive}, and \autoref{sec:acceptStop} describes the acceptance and termination criteria. A construction heuristic for generating an initial feasible solution is presented in \autoref{sec:constrHeuristic}, and a brief overview of the evaluation mechanism is given in \autoref{sec:evalution}; for details, we refer to \textcite{hansenColumnGenerationbasedFixedpoint2026}. Finally, the complete ALNS algorithm is summarized in \autoref{sec:algorithm}.

\subsection{Neighborhoods} \label{sec:neighborhoods}
In an ALNS algorithm, neighborhoods are defined by the combination of destroy and repair operators. The two operator types are typically selected independently, although they may alternatively be combined into predefined pairs \parencite[see, e.g.,][]{canca_adaptive_2017}. In this work, we adopt the former approach and propose a set of destroy and repair operators as described below. All operators act on a single transport mode $\mode \in \ModeSet$ (bus), but the approach can be extended to work with multiple modes. 

\subsubsection{Destroy operators}
\begin{list}{}{\leftmargin=0pt}
\item \textit{Remove-random-line operator}:
Let $\linePoolSol_{\mode} \subseteq \linePoolSol$ denote the set of currently
operated bus lines, excluding fixed lines. A removal fraction
$\rho \sim \mathcal{U}[0, \zeta_1]$ is sampled and a random subset
$S \subseteq \linePoolSol_{\mode}$ of size
$|S| = \max\{1, \lfloor\rho\,|\linePoolSol_{\mode}|\rfloor\}$ is removed from \linePoolSol.
The $\max\{1,\cdot\}$ ensures that at least one line is always removed.
This operator primarily serves to diversify the solution.

\item \textit{Remove-worst-line operator}:
For each line $l \in \linePoolSol_{\mode}$, we compute the utilization score $u_l=\text{load}_l/\text{cap}_l$, where $\text{load}_l$ is the total passenger flow on line $l$ and $\text{cap}_l$ is the capacity of line $l$.
We then remove $l^\star=\arg\min_{l\in\linePoolSol} u_l$, i.e., the least utilized line in the current solution. By targeting the least-utilized service, this operator focuses the search by removing lines that represent an inefficient use of resources. 

\item \textit{Remove-lines-in-area operator}:
An area $A$ is drawn uniformly at random from a predefined partition of the
service region. Let $\linePoolSol_{\mode}(A) = \{l \in \linePoolSol_{\mode} : l \text{ traverses } A\}$ be the non-fixed lines that cross $A$. A random subset $S_A \subseteq \linePoolSol(A)$ of size $|S_A| = \max\{1, \lfloor\zeta_2\,|\linePoolSol_{\mode}(A)|\rfloor\}$
is removed in full (i.e.\ the entire line, not only the portion within $A$). This operator creates structured removals in parts of the network,
thereby diversifying the search.  

\item \textit{Shorten-line operator}:
A random non-fixed line $l\in \linePoolSol_{\mode}$ with at least three stops is selected, and a terminal direction (start or end) is chosen at random. If $|E_l|$ is the
number of edges of $l$, we remove $r$ terminal edges with $r \sim \mathcal{U}\{1, \max(1, \lfloor\zeta_3(|E_l|-1)\rfloor)\}$, keeping
the remaining connected subpath. The frequency remains unchanged, while the vehicle assignment is updated to reflect the reduced round-trip time.
Compared to full removal, this operator preserves part of the existing structure and generates new line variants.

\end{list}
As mentioned, these operators focus only on lines of mode $\mode$ and do not affect fixed lines (e.g., rail infrastructure). 

\subsubsection{Repair operators}
\begin{list}{}{\leftmargin=0pt}

\item \textit{Add-random-line operator}:  
Two stops $i,j\in\stopSet$ are sampled, and the shortest path between them in $G_{\mode}$ is inserted as a line $l$ if this found line is not currently part of \linePoolSol (unique) and contains at least three stops. A frequency is drawn uniformly from the set of allowable frequencies $\frequencySet_{l}$, and the vehicle assignment is set accordingly. This operator introduces diverse lines at low computational cost. 

\item \textit{Add-lines-on-backbone-flow operator}:  
Adds up to $\zeta_4$ new lines by solving a MILP for a simplified line planning problem considering an extended line pool. 
The passenger assignment is simplified by fixing passenger routes to a \textit{backbone flow}, i.e., the shortest paths between all OD pairs
$\od \in \ODset$ ignoring capacity constraints and the current line concept. Since the backbone flow relies only on the PTN it can be precomputed. 
The extended line pool consists of the lines currently in $\linePool$ together with the one-stop terminal extensions of all lines in $\linePoolSol_{\mode}$ (prepending
or appending a single adjacent stop to each end). Fixed lines do not have variants.  
The MILP selects, for each
existing line in $\linePoolSol$, exactly one variant from its extension group (or the original
line itself) and may add up to $\zeta_4$ additional lines. 
The frequency of existing lines is held fixed; new lines receive vehicle bounds derived from the minimum and maximum allowable frequency.
The resulting frequency of a new line is implicitly determined by the vehicle assignment chosen by the MILP. Details are provided in \autoref{apx:RepairMILP}.

\item \textit{Extend-line operator}:  
A non-fixed line $l\in \linePoolSol_{\mode}$ is selected at random, and $l$ is extended to a candidate stop. The candidate stops are all stops that are direct neighbors of either terminal of $l$ in the mode-specific edge set of $G_{\mode}$ and are not already part of $l$. A random stop is selected and either prepended or appended to the stop-sequence of $l$. The frequency is preserved, and the vehicle assignment is updated for the longer round-trip.

\end{list}
When adding or extending lines, newly generated lines must be distinct from existing lines in the current solution.

\subsection{Local search} \label{sec:localsearch}
Frequency setting, i.e., determining frequencies for the selected lines, is a key part of the line planning problem as it directly affects which connections are attractive to passengers through waiting and transfer times, as well as the cost of operating a line through the required number of vehicles.
The frequency-setting problem is computationally challenging in its own right and is often studied as a standalone optimization problem \parencite[see, e.g.,][]{gallo_transit_2011, Gatt2022, de-los-santos_railway_2017}.  
In the line planning literature, frequency decisions are either integrated directly into the line planning problem, typically within exact solution approaches such as column generation \parencite{borndorfer_column-generation_2007, bertsimas_data-driven_2021}, or addressed as a separate step within iterative heuristics or matheuristic frameworks \parencite{liu_matheuristic_2020, canca_integrated_2019}.  
When passenger mode choice is also considered, the frequency-setting problem becomes even more complex. In our ALNS algorithm, we handle frequency setting through a dedicated local search operator that iteratively adjusts frequencies by ranking different alternative frequency settings based on their impact on the objective function. \autoref{pseudo:localsearch} presents the local search procedure in detail, and \autoref{pseudo:frequency_oracle} describes how the line and frequency potentials are estimated. 

The local search step iteratively adjusts line frequencies over at most $\bar{k}$ iterations.  At each iteration, candidate frequency changes are evaluated based on estimated marginal impact on objective value and combined into subsets of size at most $\bar{n}$ lines. These candidate frequency adjustments are applied and fully re-evaluated in descending order of estimated potential. Improving solutions are accepted, while non-improving candidates are discarded. The search terminates when no further improvement is possible or when the iteration limit is reached. We note that, if the subset size $\bar{n}$ equals the number of lines $|\linePoolSol|$, and the marginal-improvement pre-screening step is omitted, the procedure reduces to a full frequency enumeration, effectively corresponding to a brute-force search over all possible adjustments.

\begin{algorithm}[htbp]
\caption{Local search}
\label{pseudo:localsearch}
\KwIn{Current solution $\curSol$, objective $\fcurSol$, max iterations $\bar{k}$, max simultaneous frequency changes $\bar{n}$, Instance data $\mathcal{I}$}
\KwResult{Best solution $\bestSol$ found during search}

$\bestSol \gets \curSol,\quad \fbestSol \gets \fcurSol,\quad k \gets 0,\quad \textit{improved} \gets \text{true}$\;

\While{$k < \bar{k}$}{
    $k \gets k + 1$\;
    \tcc{Re-evaluate potentials on first iteration or after an improvement}
    \If{$k = 1$ \textbf{or} \textit{improved}}{
        $P \gets \textsc{FrequencyPotential}(\mathcal{I}, \curSol)$\;
        
        Generate all subsets $\mathcal{Q}$ of up to $\bar{n}$ line/frequency pairs $(\linePT,f)$ with $P(\linePT,f) > 0$, ordered by decreasing $\sum P$\;
    }

    \If{$\mathcal{Q}$ is empty}{
        \textbf{break}\;
    }
    \tcc{Apply and evaluate next candidate in subset}
    Select best subset $Q \gets \textit{pop}(\mathcal{Q})$\;

    $\newSol \gets \curSol$\;

    \For{each $(\linePT,f) \in Q$}{
        Set frequency of $\linePT$ to $f$ in $\newSol$\;
        Update vehicle assignment of $\linePT$\;
    }

    $\fnewSol \gets \textsc{Evaluate}(\newSol)$\;

    $\textit{improved} \gets \text{false}$\;
    \tcc{Accept if improved}
    \If{$\fnewSol < \fcurSol$}{
        $\curSol \gets \newSol,\quad \fcurSol \gets \fnewSol,\quad \textit{improved} \gets \text{true}$\;

        \If{$\fnewSol < \fbestSol$}{
            $\bestSol \gets \newSol,\quad \fbestSol \gets \fnewSol$\;
        }
    }
}

\Return $\bestSol, \fbestSol$\;
\end{algorithm}

\begin{algorithm}[htbp]
\caption{Frequency potentials}
\label{pseudo:frequency_oracle}
\KwIn{Instance data $\mathcal{I}$, solution $\curSol$}
\KwResult{Improvement potentials $P(\linePT,f)$}

Initialize $P(\linePT,f) \gets 0$ for all lines $\linePT$ and frequencies $f$\;

\For{each line $\linePT \in \linePoolSol$}{
    Let $f_0$ be current frequency of $\linePT$\;

    \For{each candidate frequency $f \in \frequencySet_\linePT \setminus \{f_0\}$}{

        Compute vehicle cost change: $\Delta C^{v}(\linePT,f)$\;

        Initialize $\Delta C^{p} \gets 0$, $\Delta W \gets 0$\;

        \For{each \od pair $\od$ affected by $\linePT$}{

            Compute reduction in generalized travel cost:
            $\Delta c_{\od}(\linePT,f)$\tcp*{waiting/transfer savings from changed frequency}

            Set updated PT utility:
            $u^{PT}_{\od} \gets u^{PT}_{\od} - \Delta c_{\od}(\linePT,f)$\;

            Compute demand change via logit model:
            $\Delta q_\od$\;

            Update $\Delta C^{p} \mathrel{+}= \Delta c_\od(\linePT,f)$\quad and\quad
            $\Delta W \mathrel{+}= \Delta q_\od$\;
        }

        $\Delta Obj(\linePT,f)
        =  \Delta C^{p}
        + \Delta C^{v}
        - (\ticketFare + \subsidy)\, \Delta W$\;

        $P(\linePT,f) \gets -\Delta Obj(\linePT,f)$\;
    }
}

\Return $P$\;
\end{algorithm}

\subsection{Adaptive mechanism} \label{sec:adaptive}
The sets of destroy and repair operators are denoted by \destroySet and \repairSet, respectively. 
Each operator is associated with a modifiable weight, denoted by
\weightDestroy{} for destroy operators and \weightRepair{} for repair operators.

At each iteration, one destroy operator \destroyOperator $\in$ \destroySet and one repair operator \repairOperator 
$\in$ \repairSet is selected using a roulette-wheel selection mechanism based on the weights. 
The probability of selecting an operator is: 
\begin{equation}
    \probDestroy{i}
    =
    \frac{\weightDestroy{i}}
    {\sum_{k \in \destroySet} \weightDestroy{k}},
    \qquad
    \probRepair{j}
    =
    \frac{\weightRepair{j}}
    {\sum_{k \in \repairSet} \weightRepair{k}}.
\end{equation}

The operator weights are updated periodically after a block of $I$ iterations in order to reflect their recent performance.
All weights are initialized to $1$.
For a destroy operator $i \in \destroySet$, the weight update rule is
\begin{equation}
    \weightDestroy{i}
    =
    (1 - \reactionFactor)\weightDestroy{i}
    +
    \reactionFactor
    \left(
    \frac{\score{i}}{\attempts{i}}
    \right),
\end{equation}
where $\reactionFactor \in [0,1]$ is the reaction factor,
\score{i} is the accumulated score of operator $i$ over the last block of $I$ iterations,
and \attempts{i} is the number of times operator $i$ was selected during this period.
The same update rule is applied to repair operators $j \in \repairSet$.

The score \score{} resets to zero after each block of iterations, and is then incrementally updated each time the specific operator is selected, following the following reward scheme. 
A reward of \reward{1} is awarded if the neighborhood solution \newSol is better than the best solution \bestSol; a score of \reward{2} is awarded if it is better than the current solution \curSol; a score of \reward{3} is awarded if it is accepted even if it is worse than the incumbent. The reward parameters satisfy $\reward{1} \ge \reward{2} \ge \reward{3}$.  

\subsection{Acceptance and termination criteria} \label{sec:acceptStop}

A candidate solution \newSol is accepted as the current solution according to the Metropolis acceptance criteria from Simulated Annealing. Specifically, a solution \newSol is accepted with probability $e^{-\frac{\fnewSol - \fcurSol}{T}}$, even if it is worse than the current solution \curSol. Here, $T$ is the temperature, which is initialized as $T_0$ and decreases as $T = \coolingFactor T$ using a cooling factor $\coolingFactor \in [0,1]$. 

The algorithm is terminated after a fixed amount of computation time. 

\subsection{Construction heuristic} \label{sec:constrHeuristic}
Although the ALNS can be initialized from an empty solution, providing a feasible initial solution generally improves convergence speed. In \autoref{sec:results}, we start the algorithm from both an initial solution based on the as-is operations and from a solution generated by a greedy construction heuristic. 
The latter is based on the idea of iteratively adding lines that contribute the most to the objective function, starting from an empty solution. 
In particular, the construction heuristic picks in each step that line from the provided line pool that has the highest potential in capturing additional passengers and evaluates the objective value at each of its frequencies to pick the best possible one.
The procedure terminates when no further line can improve the objective (keeping all already chosen lines fixed).

The pseudocode of the construction procedure is provided in \autoref{apx:Construction}.

\subsection{Evaluation step} \label{sec:evalution}
Once a line concept is found, the corresponding \cng network reduces to a directed graph of a much more manageable size. 
Given this fixed network configuration and the associated service frequencies, the evaluation problem consists of determining a passenger routing through the \cng network that respects constraints on capacity \eqref{con2}, PT service \eqref{conPTService}, and PT demand \eqref{conPTShare}.
This evaluation corresponds exactly to the SAMCF problem introduced in \textcite{hansenColumnGenerationbasedFixedpoint2026}. We therefore assess each candidate solution using the \qosheuristic proposed in that paper, which determines passenger routing and mode split for the given network. 
The heuristic is detailed in Figure~\ref{fig:evaluation-cg}. 

Given an OD matrix \ODset and the \cng network built on the current line concept $\curSol$, the heuristic iteratively alternates between (i) solving a capacitated MCF problem with exogenous OD demands and (ii) updating the OD‑specific demand based on the resulting service levels. In each iteration, the MCF problem is solved using column generation where a restricted master problem (RMP) is iteratively updated by a pricing problem that dynamically generates feasible passenger paths.
Solving the MCF yields a set of PT paths ${\Pset_{\od}}^{\prime} \subseteq \PsetODPT$ from which the routed PT flow for each OD pair \od $\in$ \ODset is obtained as
\[
q_\od = \sum_{p \in {\Pset_{\od}}^{\prime}} x_{p\od}.
\]
Let \(w_\od^{i}\) denote the demand associated with OD pair \(\od\) in iteration \(i\) of the \qosheuristic, with initialization \(w_\od^{0} = w_\od\). Demand is updated according to
\[
w^{i+1}_\od =
\begin{cases}
\kappa w^{i}_\od, 
& \text{if } q_\od = 0, \\[0.3em]
\kappa q_\od w^{i}_\od 
\;+\;
(1-\kappa) p_\od w_\od,
& \text{otherwise},
\end{cases}
\]
where \(p_\od\) denotes the proportion of travelers willing to use public transport as given by the demand function evaluated at the experienced service level given by \(x_{p\od}\). The parameter $\kappa\in (0,1)$ is a relaxation parameter controlling the update speed of the iterative scheme.
Intuitively, the update rule interpolates between the demand under inelastic assumptions (a lower bound) and the demand implied by the service-level function. Although the procedure does not have theoretical guarantees on solution quality, \textcite{hansenColumnGenerationbasedFixedpoint2026} demonstrate empirically that it yields solutions of high quality very quickly. 

\subsection{ALNS algorithm} \label{sec:algorithm}
The details of the ALNS algorithm are summarized in \autoref{pseudo:ALNS}.

\begin{algorithm}[htbp]
\caption{ALNS algorithm}
\label{pseudo:ALNS}
\KwIn{Initial solution $\initSol$, initial temperature $T_0$, time limit $\tau$, cooling rate $\coolingFactor$}
\KwResult{A feasible solution to the LPP with mode choice}

\tcc{Initialize parameters}
$T \gets T_0$, $\probDestroy{i} \gets \frac{1}{|\destroySet|} \;\forall i \in \destroySet$, $\probRepair{j} \gets \frac{1}{|\repairSet|} \;\forall j \in \repairSet$, $k \gets 0$, $\curSol \gets \initSol$, $\bestSol \gets \initSol$ \;

\While{elapsed\_time $< \tau$}{

    \tcc{Select and apply destroy/repair operators}
    Select destroy operator $i \in\destroySet$ using roulette wheel selection with \probDestroy{i}\;
    Select repair operator $j \in \repairSet$ using roulette wheel selection with \probRepair{j}\;

    $\newSol \gets \curSol$\;
    $\newSol \gets \textsc{destroy}(i, \newSol)$\;
    $\newSol \gets \textsc{repair}(j, \newSol)$\;

    
    \tcc{Local search (each move evaluates internally)}
    $\newSol, \fnewSol \gets \textsc{LocalSearch}(\newSol)$\;

    \tcc{Acceptance criterion}
    \If{$\fnewSol < \fcurSol$}{
        $\curSol \gets \newSol$\;
        $\fcurSol \gets \fnewSol$\;
        \If{$\fnewSol < \fbestSol$}{
            $\bestSol \gets \newSol$\;
            $\fbestSol \gets \fnewSol$\;
        }
    }
    \Else{
        Generate random $rnd \sim \mathcal{U}(0,1)$\;
        \If{$rnd < e^{-\frac{\fnewSol - \fcurSol}{T}}$}{
            $\curSol \gets \newSol$\;
            $\fcurSol \gets \fnewSol$\;
        }
    }

    \tcc{Update operator probabilities}
    \If{$k$ is a multiple of $I$}{
        Update the weights, \weightRepair{} and \weightDestroy{}, and reset scores \score{} and counts \attempts{}. 
    }

    $T \gets \coolingFactor T$\;

    $k \gets k + 1$\;
}

\Return \bestSol \;
\end{algorithm}

\section{Computational experiments} \label{sec:results}

This section presents the computational experiments. 
We first describe the experimental setup, including the parameter settings and computational setup. 
We then present a case study on the design of a bus network for the Danish city of Odense, demonstrating the practical applicability of our ALNS algorithm and providing insights into the resulting network design and service levels.

\subsection{The case study}
We study the municipality of Odense, which in 2025 had a population of approximately 210,000 inhabitants. The public transport system of Odense consists of a bus network, a single light rail line, and several regional and long-distance railway lines that cross through the city. 
For the optimization phase, we allow only changes to the bus network. However, the \cng network also includes light‑rail and railway connections, and the OD matrix contains demand for these modes. Consequently, passenger routing may involve train and light‑rail services, which are treated as fixed alternatives (with a fixed frequency) when evaluating candidate networks. 
\begin{figure}[htbp]   
    \centering
    \begin{tikzpicture}
        \node[anchor=south west, inner sep=0] (main) at (0,0)
            {\includegraphics[width=1.0\linewidth]{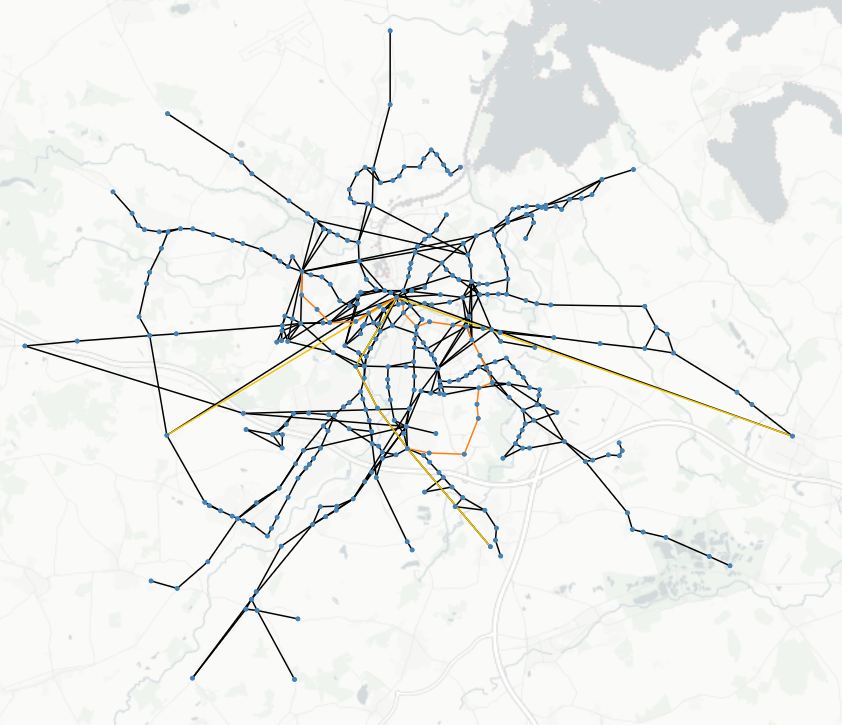}};
    \end{tikzpicture}
    \caption{The PTN network for the Odense instance. Bus edges are shown in black, train edges in yellow, and light rail edges in orange.}
    \label{fig:odense_ptn}
\end{figure}

\subsubsection{Public transport network}
The underlying network infrastructure is generated from General Transit Feed Specification (GTFS) data, which provides detailed information on stops, lines, and scheduled trips for all public transport services. 
Stop locations are extracted from the GTFS feed, and the stop sequences of individual trips are used to construct edges representing feasible, mode-specific connections between stops.  Furthermore, other relevant connections based on the road network are added manually. This initial network has 757 stops and 937 edges. 

As is common in line planning, we assume that the PTN is bidirectional. Therefore, some aggregation of the stop set is required, since the same line serves stops on both sides of the road. We use k-means to cluster stops within 100 meters into single nodes, as well as some manual adjustments to ensure meaningful groupings (e.g., at intersections or across roads) while preserving spatial separation where walking times are significant, such as between the north and south sides of the main train station. The edge set is compressed to reflect this final stop set. 
Finally, inspired by \textcite{duran-micco_designing_2022}, we simplify the network by replacing connected pairs of nodes of degree two with a single node, as any line traversing one of these nodes must also traverse the other unless one of them is a terminal node. While \textcite{duran-micco_designing_2022} compress all such nodes except terminal stops, we do not have a predefined set of terminals. Therefore, we preserve nodes where currently operated lines start or end, as well as stops with positive demand.
The resulting PTN is displayed in \autoref{fig:odense_ptn}.
The described procedure results in a PTN with 353 stops and 511 edges. The average travel time is 2.09 minutes, and the average node degree is 2.90.  

\subsubsection{Passenger costs}
Passenger costs depend on the flow through the \cng network. The arc costs are set to represent a standard generalized travel cost function, assigning different values of time to waiting, in-vehicle, and transfer arcs.
Parameter values are chosen in line with those recommended by the Danish TERESA model for socio‑economic analysis \parencite{incentive__dtu_transport_teresa_2024} and are summarized in \autoref{tab:parameterSettingPI}. Mode‑specific substitution rates for in‑vehicle time (IVT) are applied following \textcite{nielsenRelevanceDetailedTransfer2021}, with light rail classified under the local train category.

Waiting times on access arcs are modeled as half the headway of the outgoing line, consistent with the assumption of uniformly distributed passenger arrivals in the absence of timetable information. However, when headways are above a certain level, passengers tend to time their arrival to the scheduled departure time, particularly in peak hours.  
\textcite{ingvardson_passenger_2018} show that, during weekday peak hours, even with a 5-minute headway, 59\% of passengers time their arrival, while the remaining passengers arrive randomly. For headways of 10, 20, 30, and 60 minutes, the share of passengers timing their arrival is 55\%, 74\%, 90\%, and 92\%, respectively. Therefore, we adjust the waiting time calculation by weighting the random-arrival waiting time by the proportion of passengers arriving without timing their trip, while the remaining passengers are assigned a waiting time based on the hidden waiting time weight.
Transfer arc costs include both a fixed transfer penalty and the expected waiting time for the connecting service, again approximated by half the headway of the outgoing line. Egress arcs carry the ticket fare (\ticketFare), representing the monetary cost of using PT. Travel time costs are provided in \autoref{tab:parameterSettingPI}

The utility of the alternative mode is less straightforward. In practice, public transport competes with alternatives such as car, bicycle, and walking, and the relevant choice set depends on individual factors, including vehicle ownership and physical ability. We approximate the generalized cost of the alternative mode using a simplified representation.
Specifically, for each OD pair, we compute the shortest path on the PTN according to IVT. The generalized cost of the alternative mode is defined as a weighted sum of travel time and distance. We adopt the weights recommended by the TERESA model for car travel. Accordingly, travel time is scaled by a value-of-time parameter (``Car IVT''), while distance is multiplied by a unit cost per kilometer (``Car driving cost''), as reported in \autoref{tab:parameterSettingPI}.

\subsubsection{OD matrix}
\begin{figure}[htbp]
\centering

\begin{subfigure}[t]{0.48\textwidth}
\centering
\begin{tikzpicture}
    \begin{axis}[
        ybar,
        bar width=4pt,
        xlabel={Hour of day},
        ylabel={Weighted trips},
        xtick={0,3,...,21},
        ymin=0,
        enlarge x limits=0.15,
        width=0.7\linewidth,
        height=4cm,
        scale only axis,
        scaled y ticks=false,
    ]
    \addplot
    table[x=hour,y=weighted_trips,col sep=comma]
        {trips_hourly_odense.csv};
    \end{axis}
\end{tikzpicture}
\caption{Observed OD distribution}
\label{fig:obs_od_dist}
\end{subfigure}
\hfill
\begin{subfigure}[t]{0.48\textwidth}
\centering
\begin{tikzpicture}
    \begin{axis}[
        width=0.7\linewidth,
        height=4cm,
        scale only axis,
        xlabel={OD pairs},
        ylabel={Demand},
        legend style={ legend columns=2,
        at={(1,0)},
        anchor=south east},
    ]
    \addplot [line width = 1pt]  
    table[x=rank,y=cum_sum,col sep=comma]
        {cum_demand_c100.csv};


    \end{axis}
\end{tikzpicture}
\caption{Cumulative OD demand}
\label{fig:cum_od_demand}
\end{subfigure}

\caption{Characteristics of observed travel demand}
\label{fig:case_od_dist}
\end{figure}
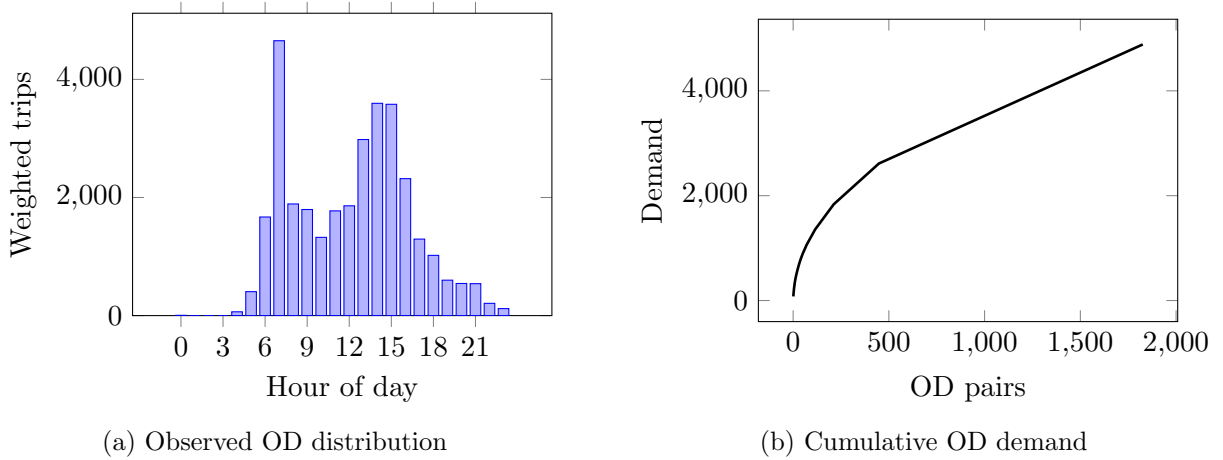

The OD matrix is generated based on data from the Danish smart card data system (Rejsekort). 
Smart card data provides precise information on travel behavior but captures only a subset of public transport trips, as some travelers use other ticket systems. To account for this, we use input from local public transport authorities to appropriately scale ridership levels. 
We construct an OD matrix for a single peak hour using data from a weekday in March 2023. Although more recent data is available, a mobile ticketing application was introduced in 2024 to complement the smart card system. As usage data from this application is not accessible, we rely on the 2023 dataset to ensure the most complete representation of demand. The OD matrix includes all trips with departure times between 7:00 and 8:00 a.m., which corresponds to the system’s peak period, as illustrated in \autoref{fig:obs_od_dist}. 
\autoref{fig:cum_od_demand} illustrates the cumulative demand distribution. 
In total, 1,826 non-zero OD pairs are identified, with an observed demand of 4,885 passenger trips.

Since this matrix reflects only observed PT demand, we calibrate a baseline total-demand matrix by inverting the logit model under current operations. 
The proportion of demand attracted to public transport is determined using the logit model in \autoref{conPTShare}. Estimating the logit parameters would require detailed behavioral data and dedicated empirical studies, and they are known to vary across OD pairs (e.g., commuters versus leisure travelers may exhibit different sensitivities).
We therefore adopt a fixed value of $\beta_{\od} = 0.05$ for all $\od \in \ODset$, which we found to provide a reasonable representation of demand sensitivity in this setting.    

For each OD pair $d$, let $q_d^{\mathrm{obs}}$ denote observed PT demand. We first compute the baseline PT share \baselineShare from \autoref{conPTShare} with \(\alpha_d = 0\), i.e., assuming no prior preference for or against PT. 
We observe that many OD pairs exhibit very low baseline probabilities, suggesting that travel decisions are influenced by factors not captured by generalized travel cost alone, such as car availability or cycling infrastructure. 
For OD pairs where the baseline mode share is less than 5\%, we set $\baselineShare = 0.05$. This lower bound avoids unrealistically small shares, which would otherwise lead to inflated total-demand estimates. 
The total demand for OD pair \od is then defined as $\widehat{w}_d = q_d^{\mathrm{obs}}/\baselineShare$.
Given $\widehat{w}_d$, we recompute the implied observed share $p_d^{\mathrm{obs}} = q_d^{\mathrm{obs}} / \widehat{w}_d$ and calibrate $\alpha_d$ such that the logit model reproduces this share under baseline conditions. Hereby, $\alpha_d$ captures OD-specific factors not represented in the generalized cost of travel.

A limitation of this approach is that it considers only OD pairs observed in the PT data and does not introduce entirely new latent OD pairs. Using this procedure, we obtain a scaled OD matrix $\hat{\ODset}$ with 59,453 
trips, representing the estimated total transportation demand for the considered period. The sensitivity of demand to service improvements is governed by $\beta_d$, meaning that improvements in service levels translate directly into changes in demand.

\autoref{apx:baseline} visualizes the baseline service probabilities under different levels of sensitivity. 

\subsubsection{Operating costs and capacities}
Operating costs are determined by the number of lines and vehicles in operation, as well as the revenue brought in from transporting passengers. As rail services are treated as fixed in our setting, we ignore the costs of these and focus on bus operations. The cost of operating a bus is based on the average hourly cost reported by the Danish Traffic Companies \parencite{nielsen_region_2024}. In addition, a fixed cost is assigned to each line. This line cost is set equal to the operating cost of one bus-hour, primarily to discourage the introduction of unnecessarily many lines.
Passenger capacity constraints are incorporated through mode-specific vehicle capacities. 
Finally, although the model includes a budget constraint on the number of vehicles and lines, this constraint is not used in the experiments. Instead, the model determines the trade-off between network investment and the service provided. 

Operating costs are offset by passenger fare revenue and public subsidies, both modeled on a per-passenger basis. Fare revenue is represented by a fixed ticket price per passenger, while the per-passenger subsidy is varied across experiments. We note that the subsidy parameter does not represent a specific reimbursement scheme, but rather controls the trade-off between ridership and operating costs in the optimization model. Selecting a realistic value for this parameter is challenging, as observed subsidy levels depend on the underlying demand and network conditions. Operational estimates suggest that the current effective subsidy translates to approximately 30--40 DKK/pax, although this value is conditional on the existing level of ridership. In the model formulation, the objective minimizes generalized travel costs both within and outside the PT system. Moreover, generalized public transport costs are modeled conservatively. For example, transfers are assumed to be uncoordinated, which can make serving some passengers less attractive than leaving them in the alternative mode. 
Consequently, without an additional per-passenger incentive, the model may favor solutions that operate only a limited public transport network and attract relatively few passengers. The subsidy, therefore, acts as a weighting parameter that encourages the model to capture demand by offsetting the cost of providing public transport service.
Given this uncertainty, we evaluate a range of subsidy levels, $[20,80]$ DKK/pax, to examine how network design and ridership outcomes respond to varying per-passenger incentives. 

An overview of all operating costs and capacity parameters is provided in \autoref{tab:parameterSettingOC}.

\begin{table}[htbp]
\centering
\begin{minipage}[t]{0.45\linewidth} 
\caption{Cost and capacity parameters for public transport operations.}
\centering
\begin{tabular}{@{}lr@{}}
\toprule
Parameter & Value\\ \midrule
Bus cost $c_{\text{bus}}$ (DKK/hour) & 880 \\
Bus capacity $\delta_{\text{bus}}$ (Pas./vehicle) & 50 \\
Light rail capacity $\delta_{\text{light rail}}$  (Pas./car) & 200 \\
Train capacity $\delta_{\text{train}}$  (Pas./car) & 300 \\
Ticket fare $\ticketFare$ (DKK/pas.) & 22 \\
\bottomrule
\end{tabular}\label{tab:parameterSettingOC} 
\end{minipage}%
\hspace{0.05\linewidth} 
\begin{minipage}[t]{0.45\linewidth} 
\caption{Parameters for generalized passenger travel costs; TERESA model  \parencite{incentive__dtu_transport_teresa_2024}}
\centering
\begin{tabular}{@{}lr@{}} 
\toprule
Parameter & Value \\ \midrule
Fixed transfer penalty (DKK) & 12 \\
Bus IVT (DKK/hour) & 119 \\
Train IVT (DKK/hour) & 107 \\
Light rail IVT (DKK/hour) & 105 \\
Waiting time (DKK/hour) & 238 \\
Hidden waiting time (DKK/hour) & 95 \\
Transferring time (DKK/hour) & 179 \\
Car IVT (DKK/hour) & 119 \\
Car driving cost (DKK/km) & 2.96 \\
\bottomrule
\end{tabular} \label{tab:parameterSettingPI}
\end{minipage}
\end{table}

\subsubsection{Lines and frequencies}
The input data includes a baseline line concept, which represents the current operations. Using the GTFS data, we identify 70 lines serving stops in the area surrounding Odense. As several of these lines differ only marginally, either by serving a small number of additional stops or diverging only outside of Odense, they are consolidated into a reduced set of 43 bidirectional, cycle-free lines. In addition, a comprehensive line pool of 797 candidate lines is generated using the tools available in the open-source LinTim program \parencite{schiewe_lintim_2024} for use within one of the ALNS neighborhoods. 
The light rail operates every 7.5 minutes during peak hours, and the train services are fixed with a headway of 20 minutes. Candidate bus headways are restricted to a discrete set $\{5, 10, 12, 15, 20, 24, 30, 40, 60\}$.

\subsection{Algorithm performance}
This section evaluates the performance of the proposed ALNS algorithm and examines the solutions produced. 
First, we examine the characteristics of the generated line plans, including changes relative to the initial solutions and the existing network. 
Second, we analyze the convergence behavior of the algorithm and assess the impact of different initialization strategies, including a chained-restart variant.  
Third, we evaluate the robustness of the obtained solutions with respect to variations in the demand elasticity parameter $\beta_{\od}$ (\autoref{sec:sensitivity}). Fourth, we investigate the contributions of individual algorithmic components through an ablation study (\autoref{sec:ablation}). Lastly, we assess the effects of alternative passenger-assignment assumptions (\autoref{sec:routing_assumptions}) and mode choice responses (\autoref{sec:modechoice}).

All experiments were performed on an 8-core Intel Xeon Gold 6226R processor with 10\,GB of memory. The ALNS algorithm was implemented in Julia version 1.12.1, and all optimization subproblems were solved using Gurobi 13.0.0.
The search procedure was terminated after three hours of computation time.

The parameter weights $\zeta_1, \zeta_2, \zeta_3,$ and $\zeta_4$ were set to $0.18$, $0.15$, $0.2$, and $2$, respectively. The reaction factor $\reactionFactor$ was set to $0.4$, and the cooling factor $\coolingFactor$ was set to $0.8$. The initial temperature was calibrated so that solutions $1\%$ worse than the initial solution are accepted with a $50\%$ probability.
Local search parameters were set to $\bar{k} = 1$ and $\bar{n} = 5$, corresponding to single-line evaluations over five iterations. The operator weight update period was set to 10 iterations. Reward parameters were fixed at $\reward{1} = 10$, $\reward{2} = 5$, and $\reward{3} = 1$.
All other parameters were tuned using the \texttt{Hyperopt.jl} package. 
The \qosheuristic is run with a maximum of 20 iterations per objective evaluation, a relaxation factor of $\kappa=0.2$, and termination criteria when the demand change per iteration falls below $0.1\%$ or the relative objective improvement falls below $0.01\%$.

\subsubsection*{Line plan and routing characteristics}
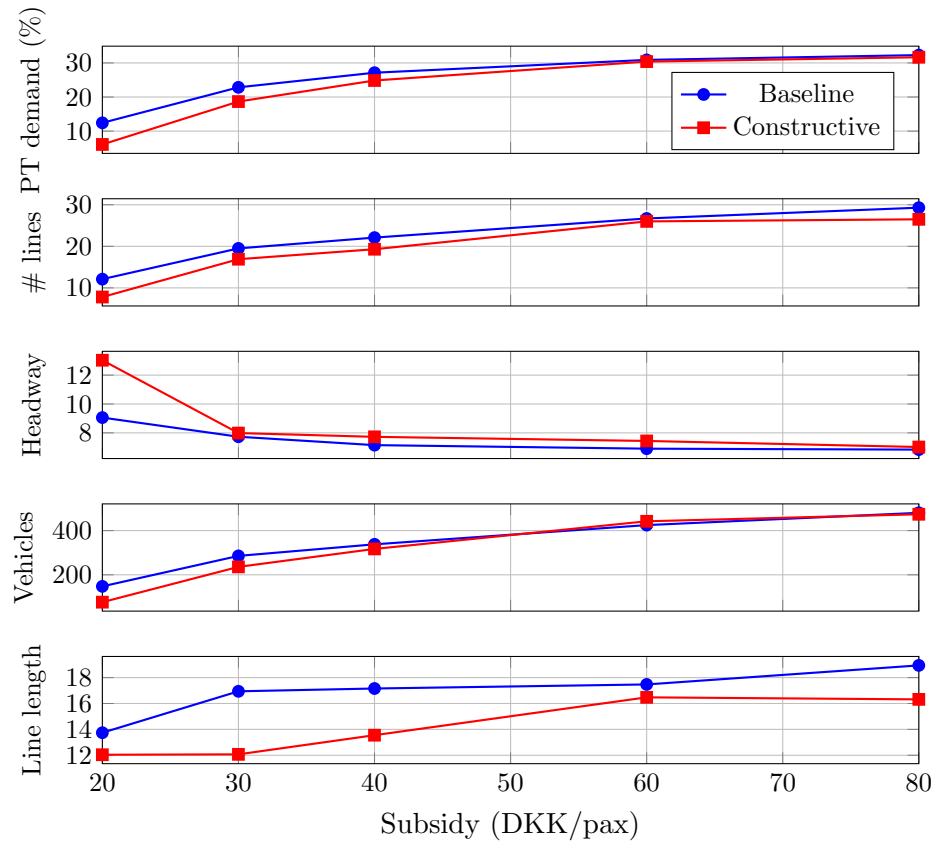
\begin{figure}[htbp]
\centering
\begin{tikzpicture}
\begin{groupplot}[
    group style={
        group size=1 by 5,
        vertical sep=0.6cm,
    },
    width=0.75\textwidth,
    height=3cm,
    grid=major,
    xlabel={Subsidy (DKK/pax)},
    ylabel style={font=\small},
    tick label style={font=\small},
    legend style={
        at={(1.02,0.5)},
        anchor=west,
        font=\small
    },
    xmin=20,
    xmax=80,
]
\nextgroupplot[
    ylabel={PT demand (\%)},
    legend pos=south east,
    xticklabels={},
    xlabel={}
]
\addplot[blue, thick, mark=*]
table[
    col sep=comma,
    x=subsidy,
    y=value_mean
]
{best_solution_panel_TOD_beta0p05_objective_combo40_pt_capture_pct_initdefault.csv};
\addlegendentry{Baseline}
\addplot[red, thick, mark=square*]
table[
    col sep=comma,
    x=subsidy,
    y=value_mean
]
{best_solution_panel_TOD_beta0p05_objective_combo40_pt_capture_pct_initconstructive.csv};
\addlegendentry{Constructive}
\nextgroupplot[
    ylabel={\# lines},
    xticklabels={},
    xlabel={}
]
\addplot[blue, thick, mark=*]
table[
    col sep=comma,
    x=subsidy,
    y=value_mean
]
{best_solution_panel_TOD_beta0p05_objective_combo40_nr_lines_initdefault.csv};
\addplot[red, thick, mark=square*]
table[
    col sep=comma,
    x=subsidy,
    y=value_mean
]
{best_solution_panel_TOD_beta0p05_objective_combo40_nr_lines_initconstructive.csv};
\nextgroupplot[
    ylabel={Headway},
    xticklabels={},
    xlabel={}
]
\addplot[blue, thick, mark=*]
table[
    col sep=comma,
    x=subsidy,
    y=value_mean
]
{best_solution_panel_TOD_beta0p05_objective_combo40_avg_headway_initdefault.csv};
\addplot[red, thick, mark=square*]
table[
    col sep=comma,
    x=subsidy,
    y=value_mean
]
{best_solution_panel_TOD_beta0p05_objective_combo40_avg_headway_initconstructive.csv};
\nextgroupplot[
    ylabel={Vehicles},
    xticklabels={},
    xlabel={}
]
\addplot[blue, thick, mark=*]
table[
    col sep=comma,
    x=subsidy,
    y=value_mean
]
{best_solution_panel_TOD_beta0p05_objective_combo40_nr_vehicles_initdefault.csv};
\addplot[red, thick, mark=square*]
table[
    col sep=comma,
    x=subsidy,
    y=value_mean
]
{best_solution_panel_TOD_beta0p05_objective_combo40_nr_vehicles_initconstructive.csv};
\nextgroupplot[
    ylabel={Line length}
]
\addplot[blue, thick, mark=*]
table[
    col sep=comma,
    x=subsidy,
    y=value_mean
]
{best_solution_panel_TOD_beta0p05_objective_combo40_avg_line_length_initdefault.csv};
\addplot[red, thick, mark=square*]
table[
    col sep=comma,
    x=subsidy,
    y=value_mean
]
{best_solution_panel_TOD_beta0p05_objective_combo40_avg_line_length_initconstructive.csv};
\end{groupplot}
\end{tikzpicture}
\caption{Characteristics of the line concepts produced by the algorithm considered different initialization strategies and subsidy levels.}
\label{fig:subsidy_lineplan_metrics}
\end{figure}

We analyze the characteristics of the resulting line plans across a subsidy range $[20,80]$ DKK/pax. 
\autoref{fig:subsidy_lineplan_metrics} presents the characteristics of the best-found line plans, averaged over 10 replications. We include the results using two different initialization strategies in the ALNS: starting from the current line concept (baseline) and starting from a solution generated by the construction heuristic (constr.).
The baseline network operates 43 lines at an average headway of approximately 41 minutes, whereas the optimized solutions include lines with an average headway between 6.8 and 13.0 minutes. 
The number of lines is reduced to between 7.8 and 29.3. This shows a clear allocation of resources in that instead of maintaining an extensive network at low frequency, the algorithm focuses on fewer, high-frequency lines. 

We note that when there is no incentive to route passengers included in the model (\subsidy = 0), the optimal solution according to the objective function is not to open any bus lines, meaning the network consists only of the fixed rail and light rail services. As the subsidy increases, additional lines are gradually introduced, but the networks are generally concentrated around fewer, longer lines than the baseline solution. 

This change in network structure is also reflected in fleet requirements. The baseline system requires 135 vehicles, whereas the optimized solutions at subsidy levels of 30 DKK or higher require larger fleets to support higher service frequencies. At lower subsidy levels, fleet sizes remain closer to the baseline, but increase steadily with the subsidy, reaching up to roughly four times the baseline fleet size. 

The demand attracted by the optimized networks is higher than that of the baseline solution, reflecting improved service quality. Initializing from the baseline solution provides the lowest objective values, with the best-found solutions capturing on average 12.4--32.3\% of total demand, depending on the subsidy level. In comparison, the current baseline system captures about 8\% of total demand (4,885 trips out of a total demand of 59,453).

For both initialization strategies, demand capture increases with the subsidy level. However, the results also indicate a clear sensitivity to the choice of initial solution. When initialized from the baseline line concept, the ALNS consistently achieves higher demand capture at low and medium subsidy levels. In contrast, initialization from the construction heuristic starts from a much sparser network and tends to keep this structure in the final solution. At higher subsidy levels (60--80 DKK), the gap between the two initialization strategies narrows, and both converge to similar levels of demand capture around 31--32\%.

This sensitivity suggests that the search process may be affected by local optima, where different starting solutions lead to distinct regions of the solution space. 
This sensitivity may be explained by the existence of many symmetric and near-equivalent solutions. Furthermore, we impose few operational constraints besides capacity constraints, which allow considerable flexibility in designing the network, increasing the size of the solution space.

\begin{figure}[htbp]
\centering
\begin{subfigure}[t]{0.49\textwidth}
\centering
\includegraphics[width=\linewidth]{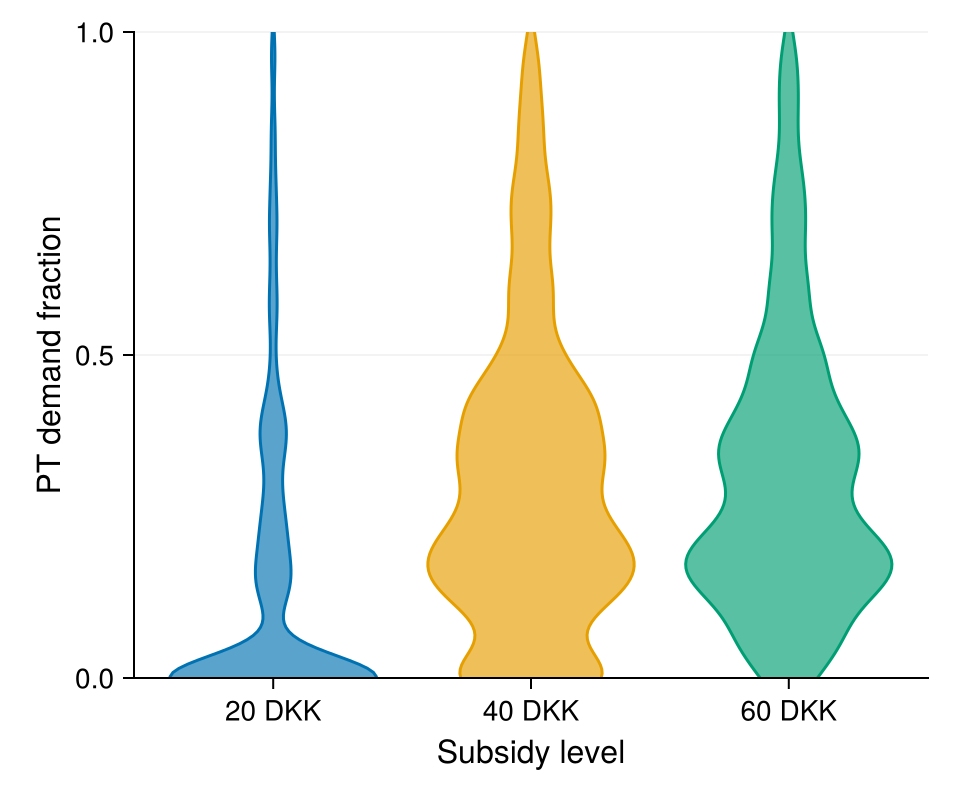}
\caption{Demand-weighted violin plot showing the distribution of PT mode share across OD pairs.}
\label{fig:pt_share_distribution}
\end{subfigure}
\hfill
\begin{subfigure}[t]{0.49\textwidth}
\centering
\includegraphics[width=\linewidth]{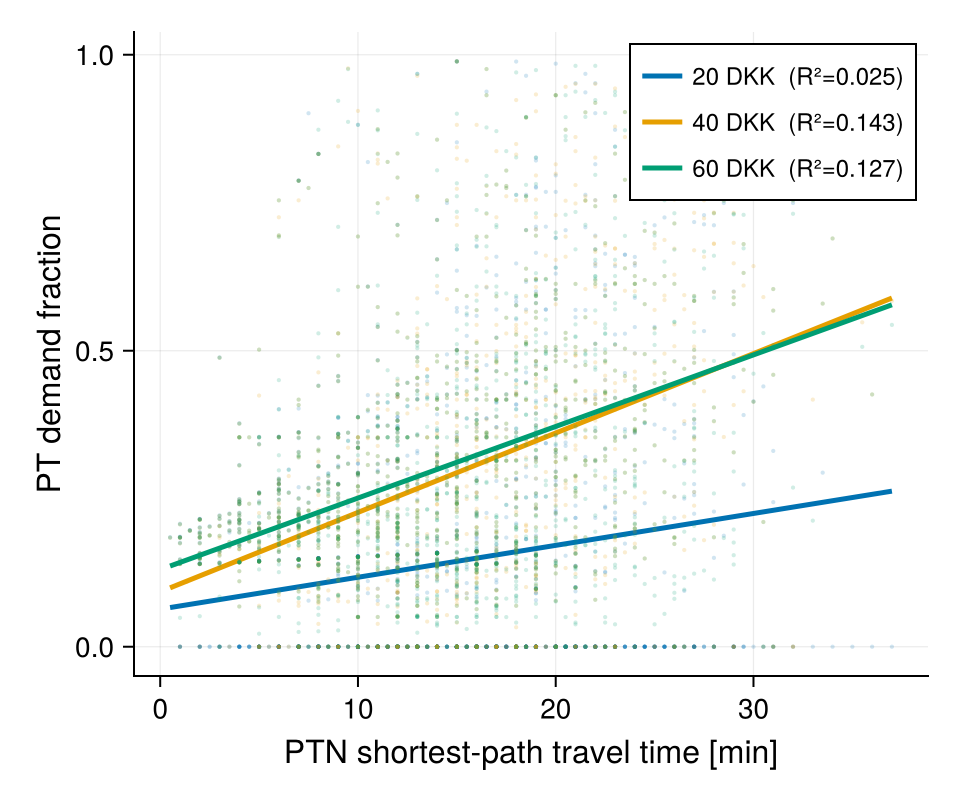}
\caption{Demand-weighted linear trend of PT fraction against PTN shortest-path travel time.}
\label{fig:pt_share_vs_distance}
\end{subfigure}
\caption{Effect of per-passenger subsidy on PT mode share. }
\label{fig:pt_share}
\end{figure}

To examine how the line concepts distribute PT service across travel demand, \autoref{fig:pt_share} analyzes the PT mode share per OD-pair, across subsidy levels 20, 40, and 60 DKK/pax.
 \autoref{fig:pt_share_distribution} shows the demand-weighted distribution of PT mode share across all OD pairs for each subsidy level, revealing whether PT service is concentrated on a few high-capture OD pairs or distributed broadly. 
In the low-subsidy setting, PT capture is highly concentrated: 60\% of OD pairs are not captured at all. The 13.6\% of passengers that are captured are distributed across the remaining OD pairs. 
The PT demand increases for the higher-subsidy settings, where the distributions for partially served OD pairs are broadly stable across subsidy levels. This suggests that increased subsidy extends coverage to previously unserved corridors rather than improving service on corridors that are already served. 
\autoref{fig:pt_share_vs_distance} further relates PT mode share to PTN shortest-path travel time. The figure shows individual OD observations together with a demand-weighted linear trend for each subsidy level.
The slope is positive across all subsidy levels, indicating that the designed networks capture a larger share of demand for longer trips. However, the relationship is weakest for low subsidy levels, where many OD pairs are not served in the PT network and therefore have a mode share of zero. 

\subsubsection*{Convergence behavior}\label{sec:convergence}
\begin{figure}[ht]
\centering
\begin{tikzpicture}


\begin{axis}[
    name=ax1,
    width=0.32\linewidth,
    height=5cm,
    xlabel={Time (s)},
    ylabel={Objective},
    title={20 DKK/pax},
    grid=major,
    scaled y ticks=false,
    yticklabel style={/pgf/number format/sci},
    xticklabel=\empty,
        ymin=2300000,
    ymax=3500000,
]

\foreach \i in {1,...,10}{
    \addplot[blue!20, thin, no markers, forget plot]
    table[col sep=comma, x=elapsed_time, y=objective_best]
{traj_sub20p0_initdefault_rep\i.csv};

    \addplot[red!20, thin, no markers, forget plot]
    table[col sep=comma, x=elapsed_time, y=objective_best]{traj_sub20p0_initconstructive_rep\i.csv};
}

\addplot[blue, thick]
table[col sep=comma, x=elapsed_time, y=objective_best_mean]
{traj_sub20p0_initdefault_mean_time.csv};

\addplot[red, thick]
table[col sep=comma, x=elapsed_time, y=objective_best_mean]
{traj_sub20p0_initconstructive_mean_time.csv};

\end{axis}

\begin{axis}[
    name=ax2,
    width=0.32\linewidth,
    height=5cm,
    at={(ax1.east)},
    anchor=west,
    title={40 DKK/pax},
    grid=major,
    xticklabel=\empty,
    yticklabels=\empty,
        ymin=2300000,
    ymax=3500000,
]

\foreach \i in {1,...,10}{
    \addplot[blue!20, thin, no markers, forget plot]
    table[col sep=comma, x=elapsed_time, y=objective_best]
{traj_sub40p0_initdefault_rep\i.csv};

    \addplot[red!20, thin, no markers, forget plot]
    table[col sep=comma, x=elapsed_time, y=objective_best]
{traj_sub40p0_initconstructive_rep\i.csv};
}

\addplot[blue, thick]
table[col sep=comma, x=elapsed_time, y=objective_best_mean]
{traj_sub40p0_initdefault_mean_time.csv};

\addplot[red, thick]
table[col sep=comma, x=elapsed_time, y=objective_best_mean]
{traj_sub40p0_initconstructive_mean_time.csv};
\end{axis}

\begin{axis}[
    name=ax3,
    width=0.32\linewidth,
    height=5cm,
    at={(ax2.east)},
    anchor=west,
    title={60 DKK/pax},
    grid=major,
    xticklabel=\empty,
    yticklabels=\empty,   
        ymin=2300000,
    ymax=3500000,
]

\foreach \i in {1,...,10}{
    \addplot[blue!20, thin, no markers, forget plot]
    table[col sep=comma, x=elapsed_time, y=objective_best]
{traj_sub60p0_initdefault_rep\i.csv};

    \addplot[red!20, thin, no markers, forget plot]
    table[col sep=comma, x=elapsed_time, y=objective_best]
{traj_sub60p0_initconstructive_rep\i.csv};
}

\addplot[blue, thick]
table[col sep=comma, x=elapsed_time, y=objective_best_mean]
{traj_sub60p0_initdefault_mean_time.csv};

\addplot[red, thick]
table[col sep=comma, x=elapsed_time, y=objective_best_mean]
{traj_sub60p0_initconstructive_mean_time.csv};

\end{axis}


\begin{axis}[
    name=ax4,
    width=0.32\linewidth,
    height=5cm,
    at={(ax1.south west)},
    yshift=-4.8cm,
    xlabel={Time (s)},
    ylabel={PT demand},
    grid=major,
    ymin=0,
    ymax=24000,
]

\foreach \i in {1,...,10}{
    \addplot[blue!20, thin, no markers, forget plot]
    table[col sep=comma, x=elapsed_time, y=pt_demand_best]
{traj_sub20p0_initdefault_rep\i.csv};

    \addplot[red!20, thin, no markers, forget plot]
    table[col sep=comma, x=elapsed_time, y=pt_demand_best]
{traj_sub20p0_initconstructive_rep\i.csv};
}

\addplot[blue, thick]
table[col sep=comma, x=elapsed_time, y=pt_demand_best_mean]
{traj_sub20p0_initdefault_mean_time.csv};

\addplot[red, thick]
table[col sep=comma, x=elapsed_time, y=pt_demand_best_mean]
{traj_sub20p0_initconstructive_mean_time.csv};

\end{axis}

\begin{axis}[
    name=ax5,
    width=0.32\linewidth,
    height=5cm,
    at={(ax4.east)},
    anchor=west,
    xlabel={Time (s)},
    grid=major,
    yticklabels=\empty,
    ymin=0,
    ymax=24000,
]

\foreach \i in {1,...,10}{
    \addplot[blue!20, thin, no markers, forget plot]
    table[col sep=comma, x=elapsed_time, y=pt_demand_best]
{traj_sub40p0_initdefault_rep\i.csv};

    \addplot[red!20, thin, no markers, forget plot]
    table[col sep=comma, x=elapsed_time, y=pt_demand_best]
{traj_sub40p0_initconstructive_rep\i.csv};
}

\addplot[blue, thick]
table[col sep=comma, x=elapsed_time, y=pt_demand_best_mean]
{traj_sub40p0_initdefault_mean_time.csv};

\addplot[red, thick]
table[col sep=comma, x=elapsed_time, y=pt_demand_best_mean]
{traj_sub40p0_initconstructive_mean_time.csv};

\end{axis}

\begin{axis}[
    name=ax6,
    width=0.32\linewidth,
    height=5cm,
    at={(ax5.east)},
    anchor=west,
    xlabel={Time (s)},
    grid=major,
    yticklabels=\empty,
        ymin=0,
    ymax=24000,
]

\foreach \i in {1,...,10}{
    \addplot[blue!20, thin, no markers, forget plot]
    table[col sep=comma, x=elapsed_time, y=pt_demand_best]
{traj_sub60p0_initdefault_rep\i.csv};

    \addplot[red!20, thin, no markers, forget plot]
    table[col sep=comma, x=elapsed_time, y=pt_demand_best]
{traj_sub60p0_initconstructive_rep\i.csv};
}

\addplot[blue, thick]
table[col sep=comma, x=elapsed_time, y=pt_demand_best_mean]
{traj_sub60p0_initdefault_mean_time.csv};

\addplot[red, thick]
table[col sep=comma, x=elapsed_time, y=pt_demand_best_mean]
{traj_sub60p0_initconstructive_mean_time.csv};

\end{axis}

\end{tikzpicture}
\caption{ALNS convergence over time starting from different initial solutions (Blue: as-is operations, Red: construction heuristic). Thin lines show the best known solution for the individual replications (10 per initialization), while thick lines show the mean trajectory.}
\label{fig:traj_all}
\end{figure}
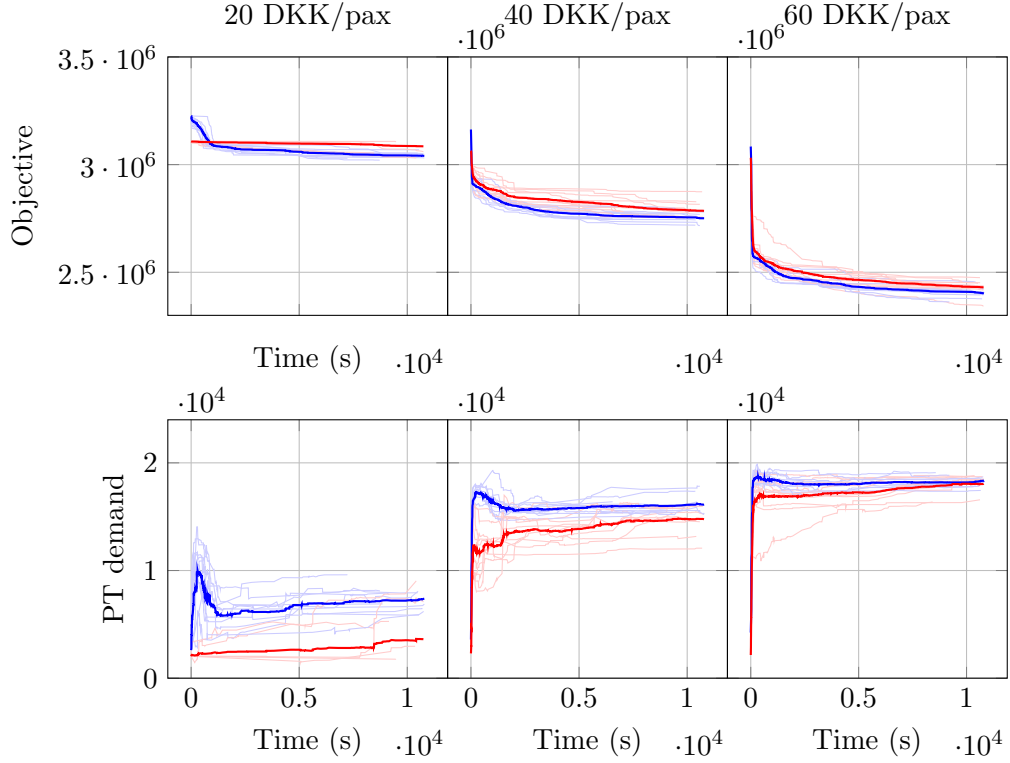

The algorithm is run for three hours, and the evolution of the best-found objective value and served PT demand is tracked over time.  \autoref{fig:traj_all} shows the trajectory of the best solution for 10 independent replications, together with the average trajectory, for different values of \subsidy. 

Over the planning horizon, the ALNS algorithm performs, on average, between approximately 1,500 and 2,800 iterations, depending on the subsidy level. 
At the lowest subsidy level, the number of iterations is considerably higher due to the sparsity of the network, which results in relatively simple routing decisions and consequently a fast evaluation step. Notably, when initialized with the construction heuristic solution, the ALNS identifies very few improving moves.
For higher subsidy values, solutions are more time-consuming to evaluate, resulting in fewer iterations explored but a higher share of improving moves (3--5\% of moves lead to a better solution). 

Most improvements occur early in the search.  Within the first 15 minutes, the objective value drops sharply while served PT demand increases correspondingly. 
After this initial phase, improvements are more gradual and occur less frequently, though they continue to accumulate throughout the full three-hour run, with improvements observed close to the time limit in most replications.

The total improvement over the run ranges from roughly 5\% at subsidy 20 DKK to over 22\% at 60 DKK for the baseline initialization, and from 1\% to 20\% for the constructive heuristic initialization. 
Despite starting from a worse objective value, the baseline initialization consistently overtakes the constructive initialization and finds better final solutions in all cases, with the final objective roughly 1--1.5\% lower on average. 
This suggests that starting from a ``better`` network constrains the search near this suboptimal area, whereas the baseline initialization explores a broader part of the solution space.

\subsubsection*{Chained restarts}
To test whether restarting from the best found solution can help the algorithm
escape local minima, we consider a chained variant with three consecutive 1-hour ALNS runs in which each run is initialized from the best solution found in the previous run. This is equivalent to periodically re-heating the temperature and reverting the incumbent to the best-known solution. We test this over 10 replications using both initialization strategies.
The overall pattern is similar across settings: most improvement occurs in the first run, while the second and third runs provide limited additional improvement.
A detailed comparison with a single 3-hour run is given in \autoref{app:chained}. Differences are negligible (never exceeding 0.52\%) and inconsistent across subsidy levels and initialization strategies. As they are comparable to the within-setting standard deviation, no systematic benefit of chained restarts can be established.

\subsection{Sensitivity analysis} \label{sec:sensitivity}
The logit parameter $\beta_{\od}$ governs how strongly passengers respond to differences in generalized travel cost between public transport and the alternative mode.
A high value of $\beta_{\od}$ implies that even modest service improvements attract substantially more passengers, while a low value reflects a slower response where mode choice is largely determined by factors other than the offered service level.
To assess the sensitivity of the algorithm and the resulting line plans to this parameter, we repeat the experiments for three values: $\beta_{\od} \in \{0.02, 0.05, 0.1\}$.
In these experiments, we consider an incentive \subsidy of 30 DKK/pax, and the algorithm is initialized from the baseline line concept, with all other settings equal to those in the other experiments.  
\begin{table}[htbp]
\centering
\caption{Line plan characteristics for solutions optimized with different mode choice parameters $\beta$. Values are reported as mean $\pm$ standard deviation over 10 replications.}
\label{tab:line_plan_characteristics}
\begin{tabular}{cccccc}
\toprule
$\beta$  & Objective & PT demand & \# lines & Avg. headway & \# vehicles \\
\midrule
$0.02$  & ${}$3,103,177 $\pm$ $2,654$ & ${}$1,460 $\pm$ $439$ & $5 \pm 1$ & $15 \pm 2$ & $35 \pm 16$ \\
$0.05$  & ${}$2,919,470 $\pm$ $18,305$ & ${}$13,578 $\pm$ $1,105$ & $20 \pm 3$ & $8 \pm 1$ & $286 \pm 32$ \\
$0.10$  & ${}$2,516,562 $\pm$ $13,289$ & ${}$30,675 $\pm$ $1,407$ & $28 \pm 3$ & $7 \pm 0$ & $472 \pm 43$ \\
\bottomrule
\end{tabular}
\end{table}
The results clearly show that the optimized line plans are highly sensitive to assumptions about passenger behavior. Compared to the baseline solution, which has an objective value of $3.20 \times 10^6$, the network design changes substantially as passenger responsiveness increases.

At low sensitivity ($\beta = 0.02$), the public transport system remains very limited, with about 5 lines, average headways of 15 minutes, and a fleet requirement of roughly 35 vehicles. In this case, only a minimal service is sustained, and both operated supply and captured demand are very low.
As $\beta$ increases to 0.05, demand rises sharply to around 13,600 passengers. This is supported by a significant expansion of the network, with the number of lines increasing to about 20 and headways decreasing to roughly 8 minutes. Fleet size also increases markedly to around 286 vehicles, reflecting the need to support higher service frequency and coverage.
At the highest value ($\beta = 0.10$), the model produces the most extensive network, with about 28 lines and demand of approximately 30,700 passengers. Service levels increase further, with headways close to 7 minutes, and fleet requirements rise to about 472 vehicles. 

Overall, relative to the baseline, the results suggest that when passengers are assumed to be more sensitive to service quality, the model responds by investing heavily in network size and capacity to reduce generalized travel costs.

\subsection{Ablation studies}\label{sec:ablation}
To assess the contribution of individual algorithmic components, we perform an ablation study in which each destroy and repair operator, as well as the local search component, is removed in turn. 
Each configuration is evaluated across 10 independent replications at subsidy levels of 20, 40, and 60 DKK/pax, starting from either initialization strategy.
\autoref{tab:ablation} reports the mean and standard deviation of the best-found objective, together with the relative difference with respect to the full ALNS configuration. 

\begin{table}[ht]
\centering
\caption{Ablation study of destroy and repair operators and the local search component. Results are aggregated across subsidy levels
and initialization strategies to assess robustness across problem settings.
\textit{Std.\ Obj.}\ is the average within-setting standard deviation across
replications, averaged over all included settings.
\textit{Diff.}\ is the percentage deviation of the mean objective from the
full configuration using all operators.}
\label{tab:ablation}
\begin{tabular}{lrrr}
\toprule
Removed component & Mean Obj. & Std. Obj. & Diff. (\%) \\ \midrule
(none) & $2.7418 \times10^{6}$ & $3.1640 \times10^{4}$ & 0.0 \\
Remove random line & $2.7415 \times10^{6}$ & $2.5680 \times10^{4}$ & -0.003 \\
Remove the least utilized line  &$ 2.7429  \times10^{6}$& $2.5953\times10^{4}$ & 0.013 \\
Extend line & $2.7567 \times10^{6}$& $2.3256\times10^{4}$ & 0.613 \\
Remove lines in area &$ 2.7595\times10^{6}$ & $2.7796\times10^{4}$ & 0.664 \\
Add backbone-flow lines & $2.7640\times10^{6} $&$ 2.8352\times10^{4}$ & 0.909 \\
Shorten line & $2.7732\times10^{6}$ & $3.0250\times10^{4} $& 1.237 \\
Local search &$ 2.7845\times10^{6}$ & $2.6648\times10^{4} $& 1.600 \\
Add random line & $2.8024\times10^{6} $& $3.1315\times10^{4}$ & 2.254\\
\bottomrule
\end{tabular}
\end{table}

Overall, most components contribute positively to solution quality, though their individual impact is modest. The random line removal operator and the worst-line removal operator have negligible effects, with the former yielding a marginally lower mean objective when removed; this may reflect overlap in functionality with the geographic destroy operator.

The largest effect is observed when removing the random line addition operator, which increases the mean objective by 2.3\%.
This suggests that the ability to introduce entirely new lines supports important diversification during the search.
The local search component and the line shortening operator have the next-largest effects, increasing the mean objective by 1.6\% and 1.2\%, respectively.
Removing the backbone-flow repair operator increases the objective by approximately 0.9\%, showing that more guided operators can benefit the search process, even when based on simplified assumptions such as fixed demand or shortest-path routing. 

Although the observed effects vary across operators, the study does not provide strong evidence for removing any individual component. Consequently, all destroy and repair operators, as well as the local search procedure, are retained in the final ALNS configuration. That said, the limited impact of some operators suggests functional overlap, indicating that further analysis could explore whether the operator set can be streamlined by reducing redundancy or by designing more complementary neighborhoods.

\subsection{Evaluation of routing assumptions}\label{sec:routing_assumptions}

The ALNS evaluates candidate line plans using a system-optimal passenger assignment, in which a central planner routes all passengers to minimize total generalized travel cost across two alternatives.
However, in practice, passengers make route choices individually, which may differ from the system optimum.
To assess the sensitivity of the results to this assumption, we re-evaluate the best-found line plans from \autoref{sec:convergence} using two alternative routing models: (i) a shortest-path assignment, in which each OD pair is assigned entirely to the path with the lowest generalized cost (SP); and (ii) a logit-based stochastic assignment, in which passengers distribute probabilistically across paths according to a logit model (MNL). We calculate the probability that a passenger chooses path $p$ from the set of available paths \PsetOD as 
\begin{equation*}
P(p) = \frac{e^{\theta c_p}}{\sum_{i \in \PsetODPT} e^{\theta c_i} + e^{\theta u^{ALT}_{\od}}} \qquad \forall p \in \PsetODPT, \od \in \ODset
\end{equation*}
where we set $\theta=-0.2$. 

\autoref{tab:samcfp_summary} summarizes the results of this analysis, showing the average generalized cost across all OD pairs under the different routing assumptions. For comparability, the evaluation is restricted to the subset of OD pairs and passengers that are (at least partially) routed in the solution. 
Overall, the results in the table show that the estimated generalized costs are quite similar between the shortest-path (SP) and system-optimal (Model) assignments. While SP consistently yields slightly lower costs, the differences are small across all subsidy levels. This indicates that, in this case, allowing for flow splitting and system-optimal assignment does not substantially increase generalized passenger costs relative to a simple shortest-path benchmark.

\begin{table}[htbp]
\centering
\caption{Average generalized PT cost for different routing assumptions. Values are averaged over 10 replications.}
\label{tab:samcfp_summary}
\begin{tabular}{lrrrr}
\toprule
\begin{tabular}[c]{@{}c@{}}Subsidy\\{[DKK]}\end{tabular}
&
\begin{tabular}[c]{@{}c@{}}Avg.\ \# routed\\passengers\end{tabular}
&
\multicolumn{3}{c}{Avg. generalized PT cost} \\
\cmidrule(lr){3-5}
& & Model [DKK] & SP [DKK] & MNL [DKK] \\
\midrule
0.0  &   560.5 & 58.49 & 58.49 & 66.58 \\
20.0 &  7386.9 & 71.07 & 71.07 & 76.95 \\
30.0 & 13578.5 & 76.24 & 76.18 & 79.36 \\
40.0 & 16129.4 & 77.75 & 77.72 & 80.09 \\
60.0 & 18363.6 & 79.30 & 79.28 & 80.13 \\
80.0 & 19201.2 & 78.79 & 78.78 & 80.82 \\
\bottomrule
\end{tabular}
\end{table}

The system-optimal assignment remains marginally higher because it distributes flows across multiple routes, some of which may involve detours or additional transfers, compared to the single minimum-cost path assumed in SP. However, the magnitude of this effect is limited in the results, suggesting that the network has sufficient capacity to route passengers on their preferred routes. 

The logit-based (MNL) assignment produces slightly higher costs than both the SP and Model assignments, but remains close in magnitude. Importantly, in the MNL evaluation, the demand is kept fixed at the level implied by the model assignment. That is, we do not model any feedback from route choice on demand or on passenger rerouting. Instead, the MNL is used purely as a post-processing step that re-evaluates generalized passenger costs conditional on existing demand levels. Extending to a fully integrated equilibrium-type model that considers mode and route choice would be an interesting next step. 

Overall, the table suggests that differences in routing assumptions have only a modest impact on estimated generalized costs in this setting, with SP providing a slightly lower bound and MNL producing higher values due to the passengers distributing across route options. 

\subsection{The impact of a mode choice response} \label{sec:modechoice}
\newcommand{\ftag}{\mathrm{fix}}
\newcommand{\dtag}{\mathrm{thres}}
\newcommand{\ltag}{\mathrm{logit}}
\newcommand{\Mdet}{\ensuremath{M_{\dtag}}\xspace}
\newcommand{\Mlogit}{\ensuremath{M_{\ltag}}\xspace}
\newcommand{\Mfix}{\ensuremath{M_{\ftag}}\xspace}
\newcommand{\Mmc}[1]{\ensuremath{M_{#1}}\xspace}      

In this section, we look at how the mode choice response, i.e., the way in which passengers choose between modes, impacts the obtained line concepts.
For this, we compare the networks produced under three alternative assumptions: a generalized travel cost threshold per OD pair (\Mdet), a variant (\Mdet(high)) of \Mdet\ which assigns a higher cost to the alternative mode, and a fixed demand response (\Mfix). 
We then re-evaluate the line concepts according to the objective in \autoref{obj}, which uses the logit function, to enable comparison. 
We refer to the model with logit-based mode choice as \Mlogit.  
The models differ only in how demand is allocated between PT and the alternative mode. 
In \Mdet, each OD pair is assigned entirely to the mode with the lower generalized cost, provided sufficient capacity is available. 
This corresponds to the limiting case of the logit model, in which the transition between the two modes becomes infinitely steep as $\beta$ increases.
As in \Mdet, the passengers also choose the mode with the lowest generalized cost in \Mdet(high), but now with an increased generalized cost of the alternative mode.
This thus represents a scenario where the competing mode is generally less attractive and helps to isolate whether differences between \Mdet and \Mlogit arise from the threshold assumption itself or from the cost parameterization. 
Model \Mfix instead approximates fixed PT demand by assigning a high penalty to unserved demand. Since the model is solved heuristically using the ALNS algorithm, the resulting line concepts may leave some demand unserved.

The implementation of the different demand responses requires only minor modifications to the evaluation procedure used as part of the ALNS algorithm. 
For \Mdet, the \qosheuristic{} is replaced by a single capacitated minimum-cost flow assignment on the \cng{} network, allowing passengers to be routed in either the PT network induced by $\linePoolSol$ or the competing mode represented by $\arcSetOptOut$. The model will assign each OD pair to the mode with the lower generalized cost, assigning passengers only to PT if there is sufficient capacity.
For \Mfix, the same minimum-cost flow assignment procedure is used, but the competing mode cost $\utilODALT$ is set to 1000 DKK/pax such that the objective function heavily penalizes unserved demand, i.e., demand assigned to the competing mode. 
To evaluate the differences in mode choice observed for these alternative demand responses, all resulting line concepts are finally re-evaluated using the \qosheuristic{} according to the objective in \autoref{obj}.

\begin{table}[] 
\scriptsize 
\centering 
\caption{Objective value, PT demand, and overall network characteristics for line concepts optimized under \Mlogit, \Mdet, \Mdet\ (high), and \Mfix considering a fixed subsidy of 40 DKK/pax. Objective value and PT demand are reported after re-evaluation under \Mlogit. Demand change reports the percentage change in PT demand between the demand assumed during optimization and the PT demand obtained after re-evaluation under \Mlogit. Objective values are given in kDKK.} \label{tab:summary_routing} 
\begin{tabular}{lrrrrrr} 
\toprule scenario & Objective & PT demand & Demand change & Lines & Vehicles & Avg. headway (min) \\ 
\midrule
\Mlogit & 2,768.0 & 15,817.1 & - & 22.1 & 339.7 & 7.3  \\ 
\Mdet & 3,052.2 & 2,285.6 & $+97.6\%$ &  6.4 & 45.1 & 14.2 \\ 
 \Mdet\ (high) & 2,816.3 & 19,131.9& $-63.9\%$ & 31.6 & 558.8 & 6.8\\
\Mfix & 2,970.8 & 20,241.2  & $-63.5\%$& 48.7 & 778.2 & 9.1 \\ 
\bottomrule 
\end{tabular} 
\end{table}

Detailed results are presented in \hyperref[app:routing]{Appendix~\ref*{app:routing}}, while \autoref{tab:summary_routing} summarizes the main findings with respect to objective value and routed demand (re-evaluated), as well as overall network characteristics. \Mdet produces very sparse networks because demand shifts to PT only when it is the cheapest option. 
When re-evaluated according to \Mlogit, the resulting networks attract more PT passengers (approximately twice as many) than under the assumed threshold response. This increase is due to the gradual probability of choosing PT introduced by the logit function. 
However, the resulting objectives remain worse than those obtained with \Mlogit.
In contrast, \Mdet(high) and \Mfix tend to produce larger networks, as they allocate resources based on demand assumptions that differ from the logit response. Both allocate resources to travelers who would be unlikely to choose PT under the logit model. 
For \Mdet(high) and \Mfix, this leads to substantial reductions in realized demand when the networks are re-evaluated using \Mlogit. 
A direct comparison of the models is not entirely straightforward, as they are based on different assumptions regarding passenger behavior and demand allocation. Nevertheless, across all settings, \Mlogit consistently achieves the best objective values. This suggests that explicitly accounting for a gradual demand response within the optimization framework can lead to more balanced network designs than approaches based on thresholds or fixed-demand assumptions.

\section{Conclusion}
\label{sec:conclusion}
In this paper, we study the line planning problem with endogenous demand, in which passenger mode choice between public transport and a competing mode is modeled using a logit function. The interaction between network design, passenger assignment, and demand estimation makes the problem substantially more challenging than the traditional line planning problem. To address this, we propose an ALNS algorithm that combines problem-specific destroy and repair operators with a local search procedure for frequency setting. A dedicated evaluation procedure is used to solve the service-aware passenger assignment problem. By generating paths and lines dynamically during the search, rather than relying on predefined candidate sets, the approach is able to scale to realistically sized networks. 

The methodology is evaluated on the regional bus network in Odense, Denmark, consisting of approximately 1,800 OD pairs. The results show that the proposed approach can generate improved line plans within a reasonable computational time for a problem size where exact optimization approaches become impractical. When endogenous demand is considered, the optimized networks differ substantially from the current operation. The solutions generally shift resources towards fewer and more frequent lines, reducing average headways from around 41 minutes in the existing network to approximately 6.8--13 minutes and increasing the captured demand share from 8\% to between 12\% and 32\% depending on the resources invested. Moreover, the resulting passenger assignment closely resembles a shortest-path routing, indicating that the optimized line plans introduce few additional detours and generally provide sufficient capacity for passengers to travel along their preferred routes.
However, the results show that the resulting network structures depend on both the initial solution and the assumed passenger response. 
In particular, the comparison of alternative demand response assumptions (specifically, a fixed and a generalized cost threshold) further demonstrates that the mode-choice assumptions influence both the resulting network design and its performance when evaluated under a common mode choice model. 
Therefore, these assumptions and their associated parameters should be considered carefully when designing and evaluating line plans.

Several limitations of the current work suggest directions for future research. Firstly, the case study could be extended in several ways. Access and egress could be represented more realistically by introducing walking connections, allowing passengers to choose between nearby stops. The line generation framework could also be extended to handle circular or looped lines. 
In the current experiments, the budget constraint is not restrictive, allowing the model to determine the preferred trade-off between service provision and ridership. Introducing tighter budget constraints would make it possible to study how much improvement can be achieved under limited resources. 
Furthermore, the framework provides a basis for evaluating policy interventions, for example, the impact of bus priority measures or fare changes. 

Secondly, there are opportunities to improve the algorithm's computational efficiency. A key computational challenge is the passenger assignment problem on the \cng network, which is solved repeatedly. Since this representation introduces a large number of auxiliary nodes, evaluating candidate solutions can become time-consuming. Using a more compact representation, such as the Direct Link Network \autocite{aktas_speeding_2024}, could improve scalability. Additional algorithmic improvements, such as more advanced destroy and repair operators or stronger diversification strategies, could also help reduce the sensitivity to the starting solution. 
Finally, developing lower-bounding techniques to assess the quality of the obtained solutions is an important direction for future research.

Thirdly, the passenger behavior model could be enhanced further. The model currently limits demand to OD pairs with observed trips. Extending the framework to account for latent demand, i.e., trips that may be induced by the improved service, would provide a more complete representation. Moreover, the model uses a single setting for the travel time sensitivity, assuming all passengers respond in the same way to changes in the service. A natural extension would be to introduce heterogeneous preferences across passenger groups, for example, different sensitivities to time, price, or frequency.

\printbibliography

@article{aktas_speeding_2024,
	title = {Speeding up the passenger assignment problem in transit planning by the {Direct Link Network} representation},
	volume = {167},
	issn = {0305-0548},
	url = {https://www.sciencedirect.com/science/article/pii/S0305054824001199},
	doi = {10.1016/j.cor.2024.106647},
	urldate = {2024-12-01},
	journal = {Computers \& Operations Research},
	author = {Aktaş, Dilay and Vermeir, Evert and Vansteenwegen, Pieter},
	month = jul,
	year = {2024},
	pages = {106647},
}

@misc{lu_line_2025,
	title = {Line planning under crowding: {A} cut-and-column generation approach},
	url = {https://arxiv.org/abs/2501.13819},
	author = {Lu, Yahan and Lieshout, Rolf N. van and Martin, Layla and Yang, Lixing},
	year = {2025},
	note = {arXiv: 2501.13819},
}

@techreport{nielsen_region_2024,
	type = {Report for the {Danish} {Regions}},
	title = {Region {Skånes} tiltag for styrkelse af den kollektive transport, og sammenligning med danske regioner ({Region} {Skåne}'s {Initiatives} to {Strengthen} {Public} {Transport}, and {Comparison} with {Danish} {Regions})},
	url = {www.regioner.dk/media/ui3p3i2x/benchmark-skaane-v16.pdf},
	urldate = {2025-05-19},
	author = {Nielsen, Otto Anker and Hylander, ens Portinson and Nyruzzaman, Robin},
	year = {2024},
}

@article{hurk_shuttle_2016,
	title = {Shuttle planning for link closures in urban public transport networks},
	volume = {50},
	issn = {0041-1655},
	doi = {10.1287/trsc.2015.0647},
	number = {3},
	journal = {Transportation Science},
	author = {Hurk, Evelien van der and Koutsopoulos, Haris N. and Wilson, Nigel and Kroon, Leo G. and Maróti, Gábor},
	month = aug,
	year = {2016},
	pages = {947--965},
}

@misc{incentive__dtu_transport_teresa_2024,
	title = {{TERESA} ({Transportministeriets} {Regnearksmodel} for {Samfundsøkonomisk} {Analyse}) for {Transportområdet} ({TERESA} – {The} {Danish} {Ministry} of {Transport}’s {Spreadsheet} {Model} for {Socio}-{Economic} {Analysis} in the {Transport} {Sector}), {Version} 6.1},
	url = {https://www.man.dtu.dk/myndighedsbetjening/teresa-og-transportoekonomiske-enhedspriser},
	urldate = {2025-06-16},
	publisher = {Transport- og Bygningsministeriet},
	author = {{Incentive \& DTU Transport}},
	year = {2024},
}

@article{bertsimas_data-driven_2021,
	title = {Data-driven transit network design at scale},
	volume = {69},
	issn = {0030-364X, 1526-5463},
	url = {https://pubsonline.informs.org/doi/10.1287/opre.2020.2057},
	doi = {10.1287/opre.2020.2057},
	language = {en},
	number = {4},
	urldate = {2024-12-02},
	journal = {Operations Research},
	author = {Bertsimas, Dimitris and Ng, Yee Sian and Yan, Julia},
	month = jul,
	year = {2021},
	pages = {1118--1133},
}

@article{ropkeAdaptiveLargeNeighborhood2006,
    title = {An adaptive large neighborhood search heuristic for the pickup and delivery problem with time windows},
    volume = {40},
    issn = {0041-1655},
    url = {https://pubsonline.informs.org/doi/abs/10.1287/trsc.1050.0135},
    doi = {10.1287/trsc.1050.0135},
    number = {4},
    urldate = {2026-04-14},
    journal = {Transportation Science},
    publisher = {INFORMS},
    author = {Ropke, Stefan and Pisinger, David},
    month = nov,
    year = {2006},
    pages = {455--472},
}

@article{ahujaSurveyVeryLargescale2002,
    title = {A survey of very large-scale neighborhood search techniques},
    volume = {123},
    copyright = {https://www.elsevier.com/tdm/userlicense/1.0/},
    issn = {0166218X},
    url = {https://linkinghub.elsevier.com/retrieve/pii/S0166218X01003389},
    doi = {10.1016/S0166-218X(01)00338-9},
    language = {en},
    number = {1-3},
    urldate = {2026-04-14},
    journal = {Discrete Applied Mathematics},
    author = {Ahuja, Ravindra K. and Ergun, Özlem and Orlin, James B. and Punnen, Abraham P.},
    month = nov,
    year = {2002},
    pages = {75--102},
}

@book{ben1985discrete,
    title = {Discrete choice analysis: {T}heory and application to travel demand},
    volume = {9},
    publisher = {MIT press},
    author = {Ben-Akiva, Moshe E and Lerman, Steven R},
    year = {1985},
}

@article{nielsenRelevanceDetailedTransfer2021,
    title = {Relevance of detailed transfer attributes in large-scale multimodal route choice models for metropolitan public transport passengers},
    volume = {147},
    issn = {09658564},
    url = {https://linkinghub.elsevier.com/retrieve/pii/S0965856421000392},
    doi = {10.1016/j.tra.2021.02.010},
    language = {en},
    urldate = {2026-05-12},
    journal = {Transportation Research Part A: Policy and Practice},
    author = {Nielsen, Otto Anker and Eltved, Morten and Anderson, Marie Karen and Prato, Carlo Giacomo},
    month = may,
    year = {2021},
    pages = {76--92},
}

@article{suman2019improvement,
	title = {Improvement in direct bus services through route planning},
	journal = {Transport Policy},
	volume = {81},
	pages = {263-274},
	year = {2019},
	issn = {0967-070X},
	doi = {https://doi.org/10.1016/j.tranpol.2019.07.001},
	url = {https://www.sciencedirect.com/science/article/pii/S0967070X19301908},
	author = {Hemant K. Suman and Nomesh B. Bolia},
}

@article{hansenExactAlgorithmPublic2026,
    title = {An exact algorithm for public transport line planning considering passenger and operational costs and lost demand},
    issn = {0377-2217},
    url = {https://www.sciencedirect.com/science/article/pii/S0377221726005886},
    doi = {10.1016/j.ejor.2026.06.034},
    urldate = {2026-07-06},
    journal = {European Journal of Operational Research},
    author = {Hansen, Siv Marie Cartland and Hoogervorst, Rowan and Nielsen, Otto Anker and Lusby, Richard Martin and van der Hurk, Evelien},
    month = jun,
    year = {2026},
}

@article{ceder_bus_1986,
	title = {Bus network design},
	volume = {20},
	doi = {10.1016/0191-2615(86)90047-0},
	number = {4},
	journal = {Transportation Research Part B},
	author = {Ceder, A. and Wilson, N.H.M.},
	year = {1986},
	pages = {331--344},
}

@article{gattermann_line_2017,
	title = {Line pool generation},
	volume = {9},
	issn = {1866-749X},
	doi = {10.1007/s12469-016-0127-x},
	number = {1-2},
	journal = {Public Transport},
	author = {Gattermann, Philine and Harbering, Jonas and Schöbel, Anita},
	month = jul,
	year = {2017},
	pages = {7--32},
}

@article{bielli2002genetic,
	title = {Genetic algorithms in bus network optimization},
	journal = {Transportation Research Part C: Emerging Technologies},
	volume = {10},
	number = {1},
	pages = {19-34},
	year = {2002},
	issn = {0968-090X},
	doi = {https://doi.org/10.1016/S0968-090X(00)00048-6},
	url = {https://www.sciencedirect.com/science/article/pii/S0968090X00000486},
	author = {Maurizio Bielli and Massimiliano Caramia and Pasquale Carotenuto}
}

@article{fan2006simulated,
	author = {Wei Fan  and Randy B. Machemehl },
	title = {Using a simulated annealing algorithm to solve the transit route network design problem},
	journal = {Journal of Transportation Engineering},
	volume = {132},
	number = {2},
	pages = {122-132},
	year = {2006},
	doi = {10.1061/(ASCE)0733-947X(2006)132:2(122)},
	URL = {https://ascelibrary.org/doi/abs/10.1061/%28ASCE%290733-947X%282006%29132%3A2%28122%29}
}

@Article{iliopoulou2022variable,
	AUTHOR = {Iliopoulou, Christina and Tassopoulos, Ioannis and Beligiannis, Grigorios},
TITLE = {A variable neighbourhood search-based algorithm for the transit route network design problem},
JOURNAL = {Applied Sciences},
VOLUME = {12},
YEAR = {2022},
NUMBER = {20},
ARTICLE-NUMBER = {10232},
URL = {https://www.mdpi.com/2076-3417/12/20/10232},
ISSN = {2076-3417},
DOI = {10.3390/app122010232}
}

@article{bussieck_optimal_1997,
	title = {Optimal lines for railway systems},
	volume = {96},
	issn = {0377-2217},
	url = {https://www.sciencedirect.com/science/article/pii/0377221795003673},
	doi = {10.1016/0377-2217(95)00367-3},
	number = {1},
	urldate = {2024-12-01},
	journal = {European Journal of Operational Research},
	author = {Bussieck, Michael R. and Kreuzer, Peter and Zimmermann, Uwe T.},
	month = jan,
	year = {1997},
	pages = {54--63},
}

@article{gatt_solving_2025,
	title = {Solving the line planning problem with service-levels using a column generation-based heuristic algorithm},
	volume = {14},
	issn = {21924376},
	url = {https://linkinghub.elsevier.com/retrieve/pii/S2192437625000135},
	doi = {10.1016/j.ejtl.2025.100164},
	language = {en},
	urldate = {2025-12-15},
	journal = {EURO Journal on Transportation and Logistics},
	author = {Gatt, Hector and Freche, Jean-Marie and Lehuédé, Fabien and Yeung, Thomas G.},
	year = {2025},
	pages = {100164},
}

@inproceedings{Gatt2022,
	address = {Dagstuhl, Germany},
	series = {{OpenAccess} series in informatics (oasics)},
	title = {A bilevel model for the frequency setting problem},
	volume = {106},
	url = {https://doi.org/10.4230/OASIcs.ATMOS.2022.5},
	doi = {10.4230/OASIcs.ATMOS.2022.5},
	booktitle = {22nd symposium on algorithmic approaches for transportation modelling, optimization, and systems ({ATMOS} 2022)},
	publisher = {Schloss Dagstuhl – Leibniz-Zentrum für Informatik},
	author = {Gatt, Hadrien and Freche, Jean-Marc and Laurent, Anne and Lehuédé, François},
	editor = {D'Emidio, Marco and Lindner, Nils},
	year = {2022},
	pages = {5:1--5:8},
}

@article{duran-micco_survey_2022,
	title = {A survey on the transit network design and frequency setting problem},
	volume = {14},
	issn = {1866-749X, 1613-7159},
	url = {https://link.springer.com/10.1007/s12469-021-00284-y},
	doi = {10.1007/s12469-021-00284-y},
	language = {en},
	number = {1},
	urldate = {2024-11-28},
	journal = {Public Transport},
	author = {Durán-Micco, Javier and Vansteenwegen, Pieter},
	month = mar,
	year = {2022},
	pages = {155--190},
}

@article{duran-micco_designing_2022,
	title = {Designing bus line plans for realistic cases - the {Utrecht} case study},
	volume = {187},
	issn = {0957-4174},
	url = {https://www.sciencedirect.com/science/article/pii/S0957417421012732},
	doi = {10.1016/j.eswa.2021.115918},
	urldate = {2024-11-28},
	journal = {Expert Systems with Applications},
	author = {Durán-Micco, Javier and van Kooten Niekerk, Marcel and Vansteenwegen, Pieter},
	month = jan,
	year = {2022},
	pages = {115918},
}

@article{liu_matheuristic_2020,
	title = {A matheuristic iterative approach for profit-oriented line planning applied to the {Chinese} high-speed railway network},
	volume = {2020},
	copyright = {Copyright © 2020 Di Liu et al.},
	issn = {2042-3195},
	url = {https://onlinelibrary.wiley.com/doi/abs/10.1155/2020/4294195},
	doi = {10.1155/2020/4294195},
	language = {en},
	number = {1},
	urldate = {2025-03-25},
	journal = {Journal of Advanced Transportation},
	author = {Liu, Di and Durán Micco, Javier and Lu, Gongyuan and Peng, Qiyuan and Ning, Jia and Vansteenwegen, Pieter},
	year = {2020},
	pages = {4294195},
}

@article{fan2008tabu,
	author = {Fan, Wei and Machemehl, Randy B.},
	title = {Tabu search strategies for the public transportation network optimizations with variable transit demand},
	journal = {Computer-Aided Civil and Infrastructure Engineering},
	volume = {23},
	number = {7},
	pages = {502-520},
	doi = {https://doi.org/10.1111/j.1467-8667.2008.00556.x},
	url = {https://onlinelibrary.wiley.com/doi/abs/10.1111/j.1467-8667.2008.00556.x},
year = {2008}
}

@article{hartleb_modeling_2023,
	title = {Modeling and solving line planning with mode choice},
	volume = {57},
	issn = {0041-1655, 1526-5447},
	url = {https://pubsonline.informs.org/doi/10.1287/trsc.2022.1171},
	doi = {10.1287/trsc.2022.1171},
	language = {en},
	number = {2},
	urldate = {2024-12-01},
	journal = {Transportation Science},
	author = {Hartleb, Johann and Schmidt, Marie and Huisman, Dennis and Friedrich, Markus},
	month = mar,
	year = {2023},
	pages = {336--350},
}

@article{goerigk_line_2017,
	title = {Line planning with user-optimal route choice},
	volume = {259},
	issn = {03772217},
	url = {https://linkinghub.elsevier.com/retrieve/pii/S0377221716308682},
	doi = {10.1016/j.ejor.2016.10.034},
	language = {en},
	number = {2},
	urldate = {2024-12-02},
	journal = {European Journal of Operational Research},
	author = {Goerigk, Marc and Schmidt, Marie},
	month = jun,
	year = {2017},
	pages = {424--436},
}

@incollection{pisinger_large_2019,
	title = {Large neighborhood search},
	volume = {272},
	issn = {08848289},
	doi = {10.1007/978-3-319-91086-4_4},
	booktitle = {Handbook of {Metaheuristics}},
	publisher = {Springer New York LLC},
	author = {Pisinger, David and Ropke, Stefan},
	year = {2019},
	pages = {99--127},
}

@article{borndorfer_column-generation_2007,
	title = {A column-generation approach to line planning in public transport},
	volume = {41},
	issn = {0041-1655, 1526-5447},
	url = {https://pubsonline.informs.org/doi/10.1287/trsc.1060.0161},
	doi = {10.1287/trsc.1060.0161},
	language = {en},
	number = {1},
	urldate = {2024-11-27},
	journal = {Transportation Science},
	author = {Borndörfer, Ralf and Grötschel, Martin and Pfetsch, Marc E.},
	month = feb,
	year = {2007},
	pages = {123--132},
}

@article{cipriani_bus_2012,
	title = {A bus network design procedure with elastic demand for large urban areas},
	volume = {4},
	issn = {1613-7159},
	url = {https://doi.org/10.1007/s12469-012-0051-7},
	doi = {10.1007/s12469-012-0051-7},
	language = {en},
	number = {1},
	urldate = {2025-07-22},
	journal = {Public Transport},
	author = {Cipriani, Ernesto and Gori, Stefano and Petrelli, Marco},
	month = jul,
	year = {2012},
	pages = {57--76},
}

@article{cancela2015mathematical,
	title = {Mathematical programming formulations for transit network design},
	journal = {Transportation Research Part B: Methodological},
	volume = {77},
	pages = {17-37},
	year = {2015},
	issn = {0191-2615},
	doi = {https://doi.org/10.1016/j.trb.2015.03.006},
	url = {https://www.sciencedirect.com/science/article/pii/S0191261515000491},
	author = {Héctor Cancela and Antonio Mauttone and María E. Urquhart}
}

@article{gallo_transit_2011,
	title = {The transit network design problem with elastic demand and internalisation of external costs: {An} application to rail frequency optimisation},
	volume = {19},
	issn = {0968-090X},
	shorttitle = {The transit network design problem with elastic demand and internalisation of external costs},
	url = {https://www.sciencedirect.com/science/article/pii/S0968090X11000349},
	doi = {10.1016/j.trc.2011.02.008},
	number = {6},
	urldate = {2025-07-22},
	journal = {Transportation Research Part C: Emerging Technologies},
	author = {Gallo, Mariano and Montella, Bruno and D’Acierno, Luca},
	month = dec,
	year = {2011},
	pages = {1276--1305},
}

@article{canca_adaptive_2017,
	title = {An adaptive neighborhood search metaheuristic for the integrated railway rapid transit network design and line planning problem},
	volume = {78},
	issn = {0305-0548},
	url = {https://www.sciencedirect.com/science/article/pii/S0305054816302039},
	doi = {10.1016/j.cor.2016.08.008},
	urldate = {2025-07-23},
	journal = {Computers \& Operations Research},
	author = {Canca, David and De-Los-Santos, Alicia and Laporte, Gilbert and Mesa, Juan A.},
	month = feb,
	year = {2017},
	pages = {1--14},
}

@article{de-los-santos_railway_2017,
	title = {The railway line frequency and size setting problem},
	volume = {9},
	issn = {1613-7159},
	url = {https://doi.org/10.1007/s12469-017-0154-2},
	doi = {10.1007/s12469-017-0154-2},
	language = {en},
	number = {1},
	urldate = {2025-07-28},
	journal = {Public Transport},
	author = {De-Los-Santos, Alicia and Laporte, Gilbert and Mesa, Juan A. and Perea, Federico},
	month = jul,
	year = {2017},
	pages = {33--53},
}

@article{canca_integrated_2019,
	title = {Integrated railway rapid transit network design and line planning problem with maximum profit},
	volume = {127},
	issn = {1366-5545},
	url = {https://www.sciencedirect.com/science/article/pii/S1366554518309979},
	doi = {10.1016/j.tre.2019.04.007},
	urldate = {2025-07-21},
	journal = {Transportation Research Part E: Logistics and Transportation Review},
	author = {Canca, David and De-Los-Santos, Alicia and Laporte, Gilbert and Mesa, Juan A.},
	month = jul,
	year = {2019},
	pages = {1--30},
}

@article{szeto_hybrid_2012,
	title = {Hybrid artificial bee colony algorithm for transit network design},
	volume = {2284},
	issn = {0361-1981},
	url = {https://doi.org/10.3141/2284-06},
	doi = {10.3141/2284-06},
	language = {EN},
	number = {1},
	urldate = {2025-07-21},
	journal = {Transportation Research Record},
	publisher = {SAGE Publications Inc},
	author = {Szeto, W. Y. and Jiang, Yu},
	month = jan,
	year = {2012},
	pages = {47--56},
}

@misc{schiewe_lintim_2024,
	title = {{LinTim}: {An} integrated environment for mathematical public transport optimization},
	url = {https://nbn-resolving.de/urn:nbn:de:hbz:386-kluedo-83557},
	publisher = {Kaiserslautern - Fachbereich Mathematik},
	author = {Schiewe, Philine and Schöbel, Anita and Jäger, Sven and Albert, Sebastian and Biedinger, Christine and Dahlheimer, Thorsten and Grafe, Vera and Roth, Sarah and Schiewe, Alexander and Spühler, Felix and Stinzendörfer, Moritz and Urban, Reena},
	year = {2024},
	note = {Pages: 156},
}

@inproceedings{schobel_line_2006,
	title = {Line planning with minimal traveling time},
	url = {http://drops.dagstuhl.de/opus/volltexte/2006/660},
	booktitle = {5th workshop on algorithmic methods and models for optimization of railways, {ATMOS} 2005},
	author = {Schöbel, Anita and Scholl, Susanne},
	year = {2006},
}

@article{schobel_line_2012,
	title = {Line planning in public transportation: {Models} and methods},
	volume = {34},
	shorttitle = {Line planning in public transportation},
	doi = {10.1007/s00291-011-0251-6},
	number = {3},
	journal = {OR Spectrum},
	author = {Schöbel, A.},
	year = {2012},
	pages = {491--510},
}

@misc{schmidt_planning_2024,
	title = {Planning and optimizing transit lines},
	url = {http://arxiv.org/abs/2405.10074},
	doi = {10.48550/arXiv.2405.10074},
	language = {en},
	urldate = {2024-11-27},
	publisher = {arXiv},
	author = {Schmidt, Marie and Schöbel, Anita},
	month = may,
	year = {2024},
	note = {arXiv:2405.10074 [math]},
}

@incollection{shawUsingConstraintProgramming1998,
    address = {Berlin, Heidelberg},
    title = {Using constraint programming and local search methods to solve vehicle routing problems},
    volume = {1520},
    url = {http://link.springer.com/10.1007/3-540-49481-2_30},
    doi = {10.1007/3-540-49481-2_30},
    urldate = {2026-04-14},
    booktitle = {Principles and practice of constraint programming — {CP98}},
    publisher = {Springer Berlin Heidelberg},
    author = {Shaw, Paul},
    editor = {Goos, Gerhard and Hartmanis, Juris and Van Leeuwen, Jan and Maher, Michael and Puget, Jean-Francois},
    year = {1998},
    note = {Series Title: Lecture Notes in Computer Science},
    pages = {417--431},
}

@article{guihaire_transit_2008,
	title = {Transit network design and scheduling: {A} global review},
	volume = {42},
	issn = {0965-8564},
	shorttitle = {Transit network design and scheduling},
	url = {https://www.sciencedirect.com/science/article/pii/S0965856408000888},
	doi = {10.1016/j.tra.2008.03.011},
	number = {10},
	urldate = {2024-11-28},
	journal = {Transportation Research Part A: Policy and Practice},
	author = {Guihaire, Valérie and Hao, Jin-Kao},
	month = dec,
	year = {2008},
	pages = {1251--1273},
}

@article{ingvardson_passenger_2018,
	title = {Passenger arrival and waiting time distributions dependent on train service frequency and station characteristics: {A} smart card data analysis},
	volume = {90},
	issn = {0968-090X},
	shorttitle = {Passenger arrival and waiting time distributions dependent on train service frequency and station characteristics},
	url = {https://www.sciencedirect.com/science/article/pii/S0968090X1830319X},
	doi = {10.1016/j.trc.2018.03.006},
	urldate = {2025-05-13},
	journal = {Transportation Research Part C: Emerging Technologies},
	author = {Ingvardson, Jesper Bláfoss and Nielsen, Otto Anker and Raveau, Sebastián and Nielsen, Bo Friis},
	month = may,
	year = {2018},
	pages = {292--306},
}

@misc{hansenColumnGenerationbasedFixedpoint2026,
    title = {A column generation-based fixed-point heuristic for the service-aware multi-commodity flow problem},
    copyright = {Creative Commons Attribution 4.0 International},
    url = {https://arxiv.org/abs/2607.10459},
    doi = {10.48550/ARXIV.2607.10459},
    urldate = {2026-07-14},
    publisher = {arXiv},
    author = {Hansen, Siv Marie Cartland and Lusby, Richard Martin},
    month = jul,
    year = {2026},
    note = {Version Number: 1},
}

\appendix
\renewcommand{\thesection}{\Alph{section}}
\section{Construction Heuristic}\label{apx:Construction}
\autoref{pseudo:constructive} presents the constructive heuristic used to generate an initial feasible solution for the line planning problem. 
The heuristic follows a greedy, iterative approach. In each iteration, it screens all candidate lines by budget feasibility, selects the line yielding the greatest incremental ridership at its highest affordable frequency, and then optimizes the frequency of the selected line considering the full objective that accounts for operating cost, ticket revenue and subsidy, and passenger generalized cost improvements. Fixed infrastructure lines (rail and light rail) are included in the solution before evaluating the candidate lines. To identify the most promising line, incremental ridership is estimated by solving a capacitated MCF problem on the CNG and subsequently bounding the demand for each OD pair according to constraints \ref{conPTShare}. Once a line has been selected, the resulting change in the objective is evaluated using the \qosheuristic.
The loop terminates when the remaining budget cannot accommodate any candidate line, when no line produces a positive ridership gain, or when the time limit is reached.

\begin{algorithm}[H]
\caption{Constructive heuristic for an initial solution}
\label{pseudo:constructive}
\KwIn{Instance data $\mathcal{I}$, budget $B$, candidate line pool $\linePool$, candidate frequencies \frequencySet, time limit $T$}
\KwResult{A feasible line concept $\initSol$}

Initialize empty solution $\initSol$, captured demand $W \gets 0$ and total operating cost $C \gets 0$\;

\If{fixed lines are required}{
    Add all fixed lines to $\initSol$ at their given frequency\; 
    Evaluate passenger demand and routing using \qosheuristic\;
    Update $C$ and $W$\;
    \If{$C > B$}{\textbf{exit} (no feasible solution exists)}
}

\While{$\linePool \neq \emptyset$ \textbf{and} elapsed time $\leq T$}{
    Set remaining budget $B_r \gets B - C$\;

    \tcc{Screen lines by budget feasibility}
    \ForEach{$l \in \linePool$}{
        \If{\text{cost of} $\linePT$ \text{at} $\min(\frequencySet_l)$ \text{exceeds} $B_r$}{
            Remove $\linePT$ from \linePool\;
        }
    }
    \If{no feasible line remains}{\textbf{break}}

    \tcc{Estimate demand impact at highest frequency within budget}
    \ForEach{remaining $l \in \linePool$}{
        Let $f^{\prime}_l$ denote the highest feasible frequency in $\frequencySet_l$ given the remaining budget $B_r$\;
        Evaluate incremental demand $\Delta W(l)$ by adding $\linePT$ at $f^{\prime}_l$\;
    }

    Select line $l^\ast \in \arg\max_{l \in \linePool} \Delta W(l)$

    \If{$\Delta W(l^\ast) \leq 0$}{\textbf{break}}

    \tcc{Optimize frequency choice for selected line}
    \ForEach{$f \in \frequencySet_{l^\ast}$}{
       \If{\text{cost of} $l^\ast$ \text{at} frequency $f$ \text{exceeds} $B_r$}{
            $\Delta Obj(l^\ast, f) \gets -\infty$}
        \Else{
            Evaluate full objective change $\Delta Obj(l^\ast, f)$ using \qosheuristic\;
        }
    }

    Select frequency $f^\ast \in \arg\max_{f \in \frequencySet_{l^\ast}} \Delta Obj(l^\ast, f)$

    \If{$\Delta Obj(l^\ast, f^\ast) > 0$}{
        Add $l^\ast$ to $\initSol$ with frequency $f^\ast$\;
        Update $C$ and $W$\;
    }

    Remove $l^\ast$ from $\linePool$\;
}

\Return $\initSol$\;
\end{algorithm}

\section{Baseline service operations}\label{apx:baseline}
This section examines the baseline service operations, which serve as a reference point for comparing the optimized line concepts.
The baseline service operations are derived from the current line concept obtained from GTFS data. 

\subsection{Demand functions for different logit-parameters}

\autoref{fig:logit_curves} illustrates the mode shares predicted by the logit model under the current line concept before calibrating the OD-specific constants ($\alpha_d$). For a fixed generalized cost of the competing mode, the logit model determines the PT mode share as a function of the generalized PT cost: a lower PT cost shifts the operating point leftward along the curve, resulting in a higher predicted PT share. The figure overlays these OD-specific curves for three values of $\beta$, with scatter points indicating the current operating point for each OD pair. For $\beta = 0.02$, the curves are relatively flat, indicating that the predicted PT share is weakly sensitive to changes in the cost ratio. As $\beta$ increases, the curves become steeper, implying that the predicted PT share becomes increasingly sensitive to changes in service quality. In the experimental results, we use a fixed value of $\beta=0.05$, but assess the sensitivity of the solutions across the range of $\beta$. 

\begin{figure}[htbp]
\centering
\begin{subfigure}[t]{0.32\textwidth}
\centering
\includegraphics[width=\linewidth]{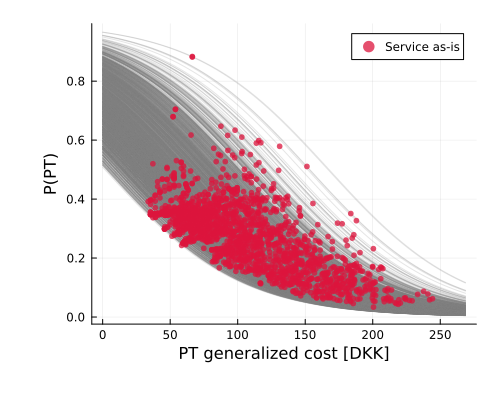}
\caption{$\beta = 0.02$}
\label{fig:logit_beta_002}
\end{subfigure}
\hfill
\begin{subfigure}[t]{0.32\textwidth}
\centering
\includegraphics[width=\linewidth]{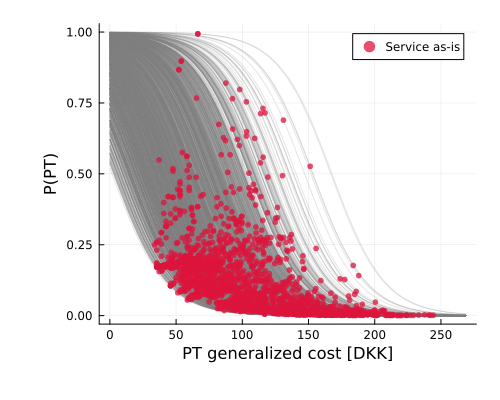}
\caption{$\beta = 0.05$}
\label{fig:logit_beta_005}
\end{subfigure}
\hfill
\begin{subfigure}[t]{0.32\textwidth}
\centering
\includegraphics[width=\linewidth]{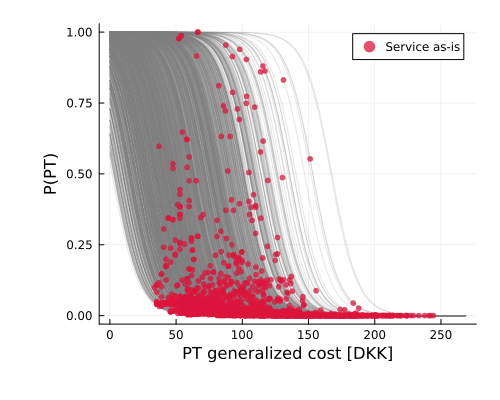}
\caption{$\beta = 0.1$}
\label{fig:logit_beta_01}
\end{subfigure}

\caption{Demand functions for each OD pair for different values of $\beta$. Scatter points indicate the current operating point for each OD pair under the baseline operations.}
\label{fig:logit_curves}
\end{figure}
\subsection{MILP-based repair operator}
\label{apx:RepairMILP}

This section details the MILP-based repair operator used in the ALNS algorithm to locally modify a line concept while preserving coverage of a given backbone passenger flow.
For each edge $e \in E$, let $d_e \ge 0$ denote the associated backbone passenger flow. We define the subset $E' := \{ e \in E \mid d_e > 0 \}$,
which can be preprocessed prior to optimization. The repair operator takes as input:
\begin{itemize}
    \item a line pool $\linePool$, where each line $l \in \linePool$ is associated with lower and upper bounds, $\underline{\phi}_l$ and $\overline{\phi}_l$, on the number of vehicles that may be assigned;
    \item the current line concept $\linePoolSol$;
    \item for each line $l \in \linePoolSol$, a set $\linePool(l) \subseteq \linePool$ of mutually exclusive line alternatives (including $l$ itself). The alternative lines $\linePool(l)$ allow for localized changes such as line extensions or line shortenings. For each alternative line $l' \in \linePool(l)$, the upper vehicle bound $\overline{\phi}_{l'}$ is set equal to the number of vehicles required to maintain the operating frequency of $l$;
    \item a target number of lines to be selected equal to $\zeta_4 + |\linePoolSol|$.
\end{itemize}
The model uses the following decision variables:
\begin{align*}
y_l &\in \{0,1\} &&\text{equals 1 if line $l \in \linePool$ is selected},\\
z_l &\ge 0       &&\text{number of vehicles assigned to line $l$},\\
u_e &\ge 0       &&\text{unserved backbone demand on edge $e \in E'$}.
\end{align*}

\begin{align}
\min \quad
& \sum_{l \in \linePool} \bigl( h_l y_l + c_{\mode(l)} z_l \bigr)
  + \sum_{e \in E'} \pi_e u_e
\label{C_obj:repair} \\[0.5em]
\text{s.t.}\quad
& \underline{\phi}_l \, y_l \le z_l \le \overline{\phi}_l \, y_l
&& \forall l \in \linePool, \label{C_con:veh_bounds} \\
& \sum_{l \in \linePool} y_l = \zeta_4 + |\linePoolSol|
&& \label{C_con:nlines} \\
& \sum_{l \in \linePool} \bigl( h_l y_l + c_{\mode(l)} z_l \bigr) \le B
&& \label{C_con:budget} \\
& \sum_{l' \in \linePool(l)} y_{l'} = 1
&& \forall l \in \linePoolSol, \label{C_con:exclusive_y} \\
& \sum_{l \in L(e)} \delta_l z_l + u_e \ge d_e
&& \forall e \in E', \label{C_con:coverage} \\
& u_e \ge 0 && \forall e \in E', \label{C_con:domain_1}\\
& z_l \ge 0 && \forall l \in \linePool,  \label{C_con:domain_2}\\
& y_l \in \{0,1\} && \forall l \in \linePool.  \label{C_con:domain_3}
\end{align}

The objective function \eqref{C_obj:repair} minimizes the total cost, comprising fixed line costs, vehicle costs, and penalties for uncovered backbone flow. 
Constraint \eqref{C_con:veh_bounds} links the line variables to the corresponding fleet assignment. If a line is selected, the number of assigned vehicles must lie between the given bounds. 
Constraint \eqref{C_con:nlines} fixes the number of selected lines to the number of lines in \linePoolSol plus the parameter $\zeta_4$.
Constraint \eqref{C_con:budget} is a budget constraint.
Constraint \eqref{C_con:exclusive_y} enforces that, for every line in \linePoolSol, exactly one line from its associated set of alternatives is selected.
Constraint \eqref{C_con:coverage} requires that the total passenger carrying capacity provided on each edge is sufficient to cover the corresponding backbone flow. Any unmet demand is captured by the non-negative slack variable $u_e$.
Finally, the remaining constraints define the domains of the decision variables.
\section{Chained ALNS: Comparison to Single-Run Results}\label{app:chained}

To assess whether sequential re-optimization (chaining) improves upon a single long run, we compare the best objective found by the chained ALNS against a single 3-hour run using the same operator configuration and initialization strategy, aggregated over 10 replications.
\autoref{tab:chained} reports the mean objective and standard deviation 
for both approaches, together with the relative improvement (positive values indicate the chained run achieves a lower objective).

\begin{table}[H]
\centering
\caption{Comparison of the chained ALNS (last round) against a single 3-hour 
run across subsidy levels and initialization strategies. 
\textit{Obj $\Delta$} is the percentage reduction in mean objective achieved 
by chaining relative to the single run; negative values indicate that the single 
run performs better. Results are averaged over 10 replications. }\label{tab:chained}
\begin{tabular}{@{}llllllr@{}}
\toprule
 &  & \multicolumn{2}{c}{Single} & \multicolumn{2}{c}{Chained} & Obj $\Delta.$\\ \cmidrule(lr){3-4}\cmidrule(lr){5-6}
subsidy & init &  Mean Obj.  &   Std. Obj.  & Mean Obj.  &   Std. Obj  &    \\
\midrule
\multirow{2}{*}{20}
& baseline
& $3.0420\times10^{6}$
& $1.7147\times10^{4}$
& $3.0462\times10^{6}$
& $1.0399\times10^{4}$
& -0.141\% \\

& constr.
& $3.0656\times10^{6}$
& $3.1259\times10^{4}$
& $3.0802\times10^{6}$
& $2.8465\times10^{4}$
& -0.477\% \\
\midrule
\multirow{2}{*}{30}
& baseline
& $2.9195\times10^{6}$
& $1.8305\times10^{4}$
& $2.9347\times10^{6}$
& $2.0334\times10^{4}$
& -0.521\% \\

& constr.
& $2.9452\times10^{6}$
& $3.8646\times10^{4}$
& $2.9343\times10^{6}$
& $2.3978\times10^{4}$
& 0.370\% \\
\midrule
\multirow{2}{*}{60}
& baseline
& $2.4031\times10^{6}$
& $2.4735\times10^{4}$
& $2.3983\times10^{6}$
& $2.0805\times10^{4}$
& 0.199\% \\

& constr.
& $2.4120\times10^{6}$
& $2.0802\times10^{4}$
& $2.4102\times10^{6}$
& $3.2169\times10^{4}$
& 0.077\% \\

\bottomrule
\end{tabular}
\end{table}

\section{Comparison of demand and mode-choice assumptions}\label{app:routing}

This section presents the results for the comparison described in \autoref{sec:modechoice}. We compare line concepts optimized under the logit mode choice model \Mlogit against those produced by the travel cost threshold model \Mdet and the fixed demand model \Mfix. The notation follows that section. 
For the experiments, we consider two values of $\subsidy$ and report results averaged over 10 replications. 
For \Mdet, we additionally report results for a variant denoted \Mdet\ (high), in which the arcs in \arcSetOptOut are set to $\utilODALT + \incentive$ for each OD pair, i.e., the true alternative mode cost plus the per-passenger revenue \incentive. 
Under \Mfix, the objective includes a penalty of 1000\,DKK for each unserved passenger. This penalty is sufficiently large that the ALNS constructs line concepts serving almost all demand, thereby approximating a fixed-demand setting. Nevertheless, because unserved demand remains feasible at a penalty, \Mfix\ does not represent a strictly fixed-demand model. 

\autoref{tab:modeling_assumption} quantifies the gap between the as-optimized and logit-repaired objectives. For \Mdet\ (high), the as-optimized objective is highly optimistic, and, when re-evaluating with the logit function, the objectives increase by 24.9\% and 87.3\% at subsidy levels 20 and 40, while repaired PT demand drops by approximately $63\%$ in both cases. 
This reflects the large share of passengers the model assumes will travel via PT that would not choose it under the logit specification. 
In contrast, \Mdet routes very little demand because the generalized cost of the alternative mode is lower than the PT cost for most OD pairs, making PT unattractive in this scenario. 
When re-evaluating with the logit function, PT demand increases substantially because the logit model assigns a non-zero probability of choosing PT, allowing some demand to still be routed.  
For \Mfix, the optimization produces line concepts that serve almost all demand in the as-optimized solution. However, once these networks are re-evaluated under \Mlogit, repaired PT demand decreases by about 63\%  and the repaired objective by 45--69\%.

\begin{table}[]
\centering
\caption{As-optimized and logit-repaired objectives and PT demand for solutions optimized under \Mdet, \Mdet\ (high), and \Mfix. ``Repaired'' denotes re-evaluation of the line concept under \Mlogit. The objective values are given in kDKK.}
\label{tab:modeling_assumption}
\resizebox{1.0\textwidth}{!}{
\begin{tabular}{@{}llrrrrrr@{}}
\toprule
subsidy & scenario & Reported obj. & Repaired obj. & Obj. $\Delta$ & Reported PT & Repaired PT & PT $\Delta$ \\
\midrule

\multirow{3}{*}{20}
& \Mdet
& $3122.3$
& $3112.1$
& $-0.3\%$
& $387.5$
& $1195.5$
& $208.5\%$ \\

& \Mdet\ (high)
& $2486.5$
& $3105.4$
& $24.9\%$
& $43924.8$
& $16132.8$
& $-63.3\%$ \\

& \Mfix 
& 10928.0
& 3413.6
& $-68.8\%$  
& 50939.6 
& 18860.5
& $-63.0\%$  \\

\midrule

\multirow{3}{*}{40}
& \Mdet
& $3079.4$
& $3052.2$
& $-0.9\%$
& $1156.4$
& $2285.6$
& $97.6\%$ \\

& \Mdet\ (high)
& $1503.3$
& $2816.3$
& $87.3\%$
& $53066.5$
& $19131.9$
& $-63.9\%$ \\

& \Mfix 
& 5402.6 
& 2970.8 
& $-45.0\% $
& 55516.1 
& 20241.2 
& $-63.5\%$ \\   

\bottomrule
\end{tabular}}
\end{table}

\autoref{tab:routing_regime_comparison} compares all four scenarios under the logit function evaluation. 
\Mlogit achieves the best objective for both settings. 
\Mdet results in the lowest operating costs, but also the lowest PT demand, leading to objective values that are 2.5\% and 10.3\% higher than \Mlogit. In contrast, both \Mdet(high) and \Mfix attract more PT demand than \Mlogit, but at higher operating costs. Of these, \Mdet(high) achieves objective values closest to \Mlogit, whereas \Mfix incurs the largest cost increases and objective values 7.3--12.4\% above \Mlogit.

\begin{table}[]
\centering
\caption{Objectives and PT demand for line concepts optimized under \Mlogit, \Mdet, \Mdet\ (high), and \Mfix all re-evaluated under \Mlogit. $\Delta$ vs \Mlogit reports the percentage difference relative to \Mlogit. The objective value and operating costs are given in kDKK.}
\label{tab:routing_regime_comparison}
\resizebox{1.0\textwidth}{!}{
\begin{tabular}{@{}llrrrrrr@{}}
\toprule
subsidy & scenario & Objective & $\Delta$ vs \Mlogit & PT demand & $\Delta$ vs \Mlogit & Oper. costs & $\Delta$ vs \Mlogit \\
\midrule

\multirow{4}{*}{20}
& \Mlogit
& $3037.5$
& -
& $8188.0$
& -
& $138.0$
& - \\

& \Mdet
& $3112.1$
& $2.5\%$
& $1195.5$
& $-85.4\%$
& $5.2$
& $-96.2\%$ \\

& \Mdet\ (high)
& $3105.4$
& $2.2\%$
& $16132.8$
& $97.0\%$
& $393.3$
& $185.0\%$ \\

& \Mfix 
& 3413.6 
& 12.4\%  
& 18860.5 
& 130.3\%  
& 760.8 
& 451.3\%  \\

\midrule

\multirow{4}{*}{40}
& \Mlogit
& $2768.0$
& -
& $15817.1$
& -
& $298.1$
& - \\

& \Mdet
& $3052.2$
& $10.3\%$
& $2285.6$
& $-85.6\%$
& $25.1$
& $-91.6\%$ \\

& \Mdet\ (high)
& $2816.3$
& $1.7\%$
& $19131.9$
& $21.0\%$
& $499.3$
& $67.5\%$ \\

& \Mfix 
& 2970.8 
& 7.3\%
& 20241.2 
& 28.0\% 
& 707.4 
& 137.3\%  \\

\bottomrule
\end{tabular}}
\end{table}

\autoref{tab:network_characteristics} shows the corresponding network structures. \Mdet produces the sparsest networks, with only 4.6--6.4 lines, as only a small share of OD pairs find PT competitive given the alternative mode. 
At the other end, \Mfix builds the largest networks by far, with 48.7--52.7 lines and fleets of 778--835 vehicles, reflecting the strong incentive to serve nearly all demand. \Mdet\ (high) also produces extensive networks, with 25.7 and 31.6 lines and fleet sizes of 444 and 559 vehicles.

\begin{table}[]
\centering
\caption{Network characteristics of line concepts optimized under \Mlogit, \Mdet, \Mdet\ (high), and \Mfix.}
\label{tab:network_characteristics}
\begin{tabular}{@{}llrrr@{}}
\toprule
subsidy & scenario & Lines & Vehicles & Avg. headway (min) \\
\midrule

\multirow{4}{*}{20}
& \Mlogit
& $13.3$
& $166.5$
& $8.8$ \\

& \Mdet
& $4.6$
& $24.3$
& $16.2$ \\

& \Mdet\ (high)
& $25.7$
& $444.2$
& $6.9$ \\

& \Mfix 
&52.7 
&834.8
&9.3\\

\midrule

\multirow{4}{*}{40}
& \Mlogit
& $22.1$
& $339.7$
& $7.3$ \\

& \Mdet
& $6.4$
& $45.1$
& $14.2$ \\

& \Mdet\ (high)
& $31.6$
& $558.8$
& $6.8$ \\

& \Mfix 
& 48.7 
& 778.2
& 9.1 \\

\bottomrule
\end{tabular}
\end{table}

For \Mfix we additionally consider the unscaled OD matrix of 4,885 observed trips. We note that \Mlogit does not open any bus lines for the unscaled OD matrix, so a direct \Mlogit baseline is not applicable. 
Relative to the setting with the scaled OD matrix, the networks generated with \Mfix open a similar number of lines but at substantially higher headways (19--21 minutes), reflecting that reduced demand lowers the need for capacity. 
\newpage

\end{document}